\documentclass[3p,final]{elsarticle}

\usepackage{amssymb}
\usepackage{amsmath}

\usepackage{dsfont}
\usepackage{hyperref}
\usepackage{cleveref}
\usepackage{amsthm}
\usepackage{numprint}
\usepackage{babel}[english]
\usepackage{caption}
\usepackage{subcaption}
\usepackage{ gensymb }
\usepackage{lineno}
\modulolinenumbers[5]
\usepackage{ragged2e}
\usepackage[dvipsnames]{xcolor}
\usepackage{natbib}
\usepackage{graphicx}
\usepackage{multicol}
\usepackage{multirow}
\usepackage{relsize}
\usepackage{derivative}
\usepackage{booktabs}
\usepackage{bm}
\usepackage{soul}  

\usepackage{cellspace} 
\begin{document}

\begin{frontmatter}



\title{Continuation strategies to mitigate convergence to low-performing local optima in topology optimization of sound transmission loss}


\author[Diepenbeek,FlandersMake]{Tom De Weer}
\ead{tom.deweer at kuleuven.be}
\author[KUL,FlandersMake]{Vanessa Cool}
\author[Diepenbeek,FlandersMake]{Elke Deckers}

\affiliation[Diepenbeek]{organization={KU Leuven Campus Diepenbeek, Department of Mechanical Engineering},
            addressline={Wetenschapspark 27}, 
            city={Diepenbeek},
            postcode={3590}, 
            country={Belgium}}

\affiliation[KUL]{organization={ KU Leuven, Department of Mechanical Engineering},
            addressline={Celestijnenlaan 300}, 
            city={Heverlee},
            postcode={3001}, 
            country={Belgium}}

\affiliation[FlandersMake]{organization={Flanders Make@KU Leuven},
            country={Belgium}}

\begin{abstract}
Dynamic topology optimization problems often suffer from convergence to low-performing local optima. This typically results in stiff designs that do not exploit dynamical phenomena such as antiresonance and decoupling. To obtain better designs, researchers often repeat their optimizations with different initial guesses. However, such reruns are computationally expensive and the required number is unknown. To quantify this problem, random initial guesses are sampled and tested for different frequencies on two case studies: 1) dynamic compliance minimization of a reinforced cantilever, which exhibits poor optima for driving frequencies below the first natural frequency, and 2) sound transmission loss maximization of a sandwich panel, which additionally sees a strong tendency towards low-performing optima at high frequencies. To address this issue, the study first divides techniques to reduce the needed number of reruns into four categories: global optimization, exclusion, relaxation and frequency‑shift methods. For the latter three, continuation strategies are proposed, illustrated, evaluated and compared on the sound transmission loss case, using Monte Carlo sampling to estimate success rates. All strategies show measurable benefits and trade‑offs. To support broader applicability, the study concludes with practical guidelines for dealing with convergence to poor local optima in dynamic topology optimization.
\end{abstract}




\end{frontmatter}

\section{Introduction}
\label{sec:introduction}

Mechanical structures in numerous engineering applications are continuously subjected to dynamic loading. \ 
This complicates their design, since the dynamic load must be incorporated in the design process.\ 
However, adapting the design to shift its eigenfrequencies away from the loading frequencies can lead to an unwanted loss of stiffness, which poses a difficult trade-off.\ 
To strike a balance in this design challenge, the seminal introduction of topology optimization~\cite{bendsoe_kikuchi_1988} in 1988 was quickly followed by its application to dynamic problems by Díaaz and Kikuchi~\cite{Diaz_Kikuchi_1992} in 1992.\ 
Since then, the subfield of dynamic topology optimization has seen research towards, e.g, numerical sensitivity issues and mode switching~\cite{lewis_overton_book1996}, large scale computations and model order reduction~\cite{Seok_Large_Scale_TO,Xiao2024,LI2021114149} and applications such as noise-reducing sandwich structures~\cite{cool_TO_VA}, resonant structures~\cite{Guzman2023_antiresonance_matching} and Micro Electro Mechanical Systems (MEMS)~\cite{GIANNINI2022104352}.\\ 

\begin{figure}[h]
    \centering
    \includegraphics[width=0.5\linewidth]{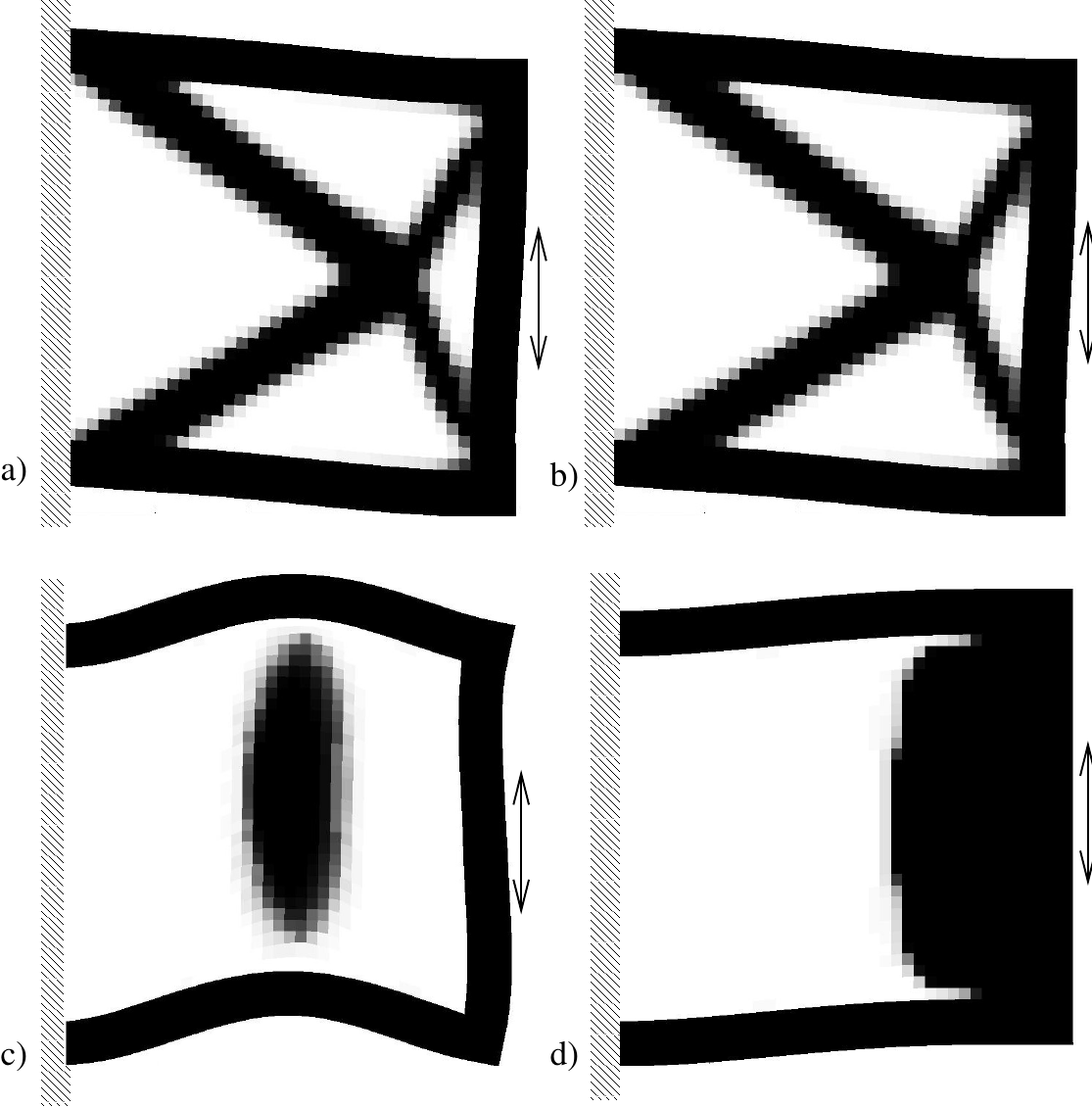}
    \caption{Figure 2.4 from Bendsøe and Sigmund~\cite{BendsøeSigmundBook2003}, used with permission from the authors and the publisher, showing four local optima of dynamic compliance minimization of a cantilever with fixed outer frame, for increasing driving frequencies.}
    \label{fig:Fig2.4_BendsoeSigmund}
\end{figure}

A particular issue in topology optimization is its many local optima.\ 
Although classical compliance minimization is already known for its strong nonconvexity, dynamic compliance minimization exacerbates this issue due to the eigenfrequencies of the structure, since they create a design space with sharp peaks between which dynamic compliance optimizers get stuck~\cite{Silva2019}.\
This is a well-known issue, as exemplified by \Cref{fig:Fig2.4_BendsoeSigmund}, which shows four local optima from Bendsøe and Sigmund~\cite{BendsøeSigmundBook2003}.\
The top designs are obtained with low driving frequencies, resulting in statically stiff structures with a high first natural frequency.\
The bottom designs are found at higher frequencies.\ 
Because they exploit an antiresonance, better dynamic performance is achieved.\

To deal with the issue of local optima in dynamic topology optimization problems, many techniques have been investigated in the literature.\
Yet, a clear overview and quantifiable comparison is lacking.\ 
This work proposes to divide the strategies into four categories.\\

The first category consists of \emph{global optimization methods}.\
Such methods are often applicable to generic optimization problems, of which design optimization is but a subcategory.\ 
A standard way is to start the optimization from a large amount of initial guesses generated randomly by, e.g., Latin hypercube sampling, and selecting the best local optimum.\
One example that improves upon this sampling concept is recent work by Herrman et al.~\cite{Herrman2024}, who use a neural network material discretization to improve the chance at a better optimum for acoustic problems.\
Another instance of this category is by Assimi and Jamali~\cite{Assimi2018}, who show that a combination of Nelder-Mead and genetic programming can lead to lighter optimal truss layouts when considering a large amount of static and dynamic loads.\
Similarly, Wang and Kai~\cite{WANG20053749} apply genetic programming on structural topology optimization problems.\
Zhang et al.~\cite{Zhang2018} apply the tunneling method to compliance minimization, moving from one optimum to a better one by iteratively adapting the objective space.\
Both Papadopoulos et al.~\cite{Papadopoulos2021} and Baeck et al.~\cite{Baeck2025} compute multiple solutions to the same topology optimization problem through a deflated barrier method that removes previously identified solutions from the search.\
Finally, recent work by Dalklint et al.~\cite{dalklint2024performanceboundstopologyoptimization} presents a method to estimate the distance from the global optimum.\ 
This provides a way to gauge performance bounds and stop global optimization techniques when they are (close to) globally optimal.\
In general, however, global optimization requires excessive computational resources, especially so in topology optimization contexts, and are therefore not further investigated in this work.\

The second category are called \emph{exclusion methods}.\ 
If the class of numerically attractive but undesired local optima is well known, constraints can be included to remove those optima from the feasible design space.\ 
Arguably the simplest example of this is a volume constraint in dynamic compliance minimization below the first eigenfrequency: without it, routines converge to a fully solid, mass-driven design.\ 
Nonetheless, in a dynamic context, exploiting dynamic behavior such as resonances by reducing mass and introducing compliance can lead to superior performance.\
In this context, a volume constraint is a numerical way to induce such behavior.\ 
A second example of an exclusion method are connectivity constraints.\ 
As highlighted in the review by Cool et al.~\cite{cool_connectivity_review}, topology optimization routines in dynamic applications often produce disconnected structures, i.e., structures with floating solid (or void) regions.\ 
Examples include broadband filter design by Dilgen and Aage~\cite{DILGEN2024104123}, band gap maximization in Mindlin plates by Halkjær et al.~\cite{Halkjaer2006} and sound transmission loss maximization of sandwich structures by Cool et al.~\cite{cool_TO_VA}.\ 
Such disconnected designs can provide superior performance but they are often undesired, for either stiffness or manufacturability reasons.\ 
Connectivity constraints thus exclude disconnected designs from the feasible design space.\
A third example are reliability methods that aim to improve the robustness of designs against stochastic variations.\
Recent numerical experiments show that such methods, by excluding non-robust designs, can improve performance even in the purely deterministic case~\cite{Elbek2025}.\ 
Generally, exclusion methods can completely prevent the appearance of undesired classes of local optima, provided that the class is understood well enough to come up with a suitable exclusion constraint.\

The third category consists of \emph{relaxation methods}.\ 
At its core, relaxation methods aim to solve a simpler (i.e., less constrained) but related problem first, and use the solution as a starting point for the actual problem.\ 
In this sense, it is strongly related to the various continuation schemes encountered in TO.\ 
For example, structural compliance minimization with the SIMP law gradually increases the penalization to promote a black-white solution.\ 
Penalization has a long history and stems back to early studies on density-based compliance minimization.\
Arguably the earliest work is that of Allaire and Kohn in 1993~\cite{allaire_kohn_1993}, which employs a grayness penalization term in the objective that increases during the optimization.\ 
The design space is ``relaxed'' by allowing density values between $0$ and $1$ and slowly transforming it to the original black-white problem by increasing the penalization.\ 
This continuation is a general feature of relaxation methods.\
More recently, Sanders et al.~\cite{Clay2020} used an objective consisting of an additional penalty term that encapsulates the error in the constitutive relations.\
The relaxation consists of the fact that the constitutive relations need not be satisfied at the optimization start.\ 
The optimization of the direct frequency response function involves choosing displacements, stresses and inertial forces on top of the usual density field.\  
It thereby avoids optimization of a functional solely of the displacements, which yields superior designs.\
Guo et al.~\cite{Guo2001} employ a second order smooth-extended technique in combination with the epsilon-relaxed method to allow for better truss layouts when considering local stress and buckling constraints.\
In general, relaxation is a powerful tool in optimization. \ 
Many techniques and variations exist and each one requires careful tuning to ensure the initial relaxation is gradually restricted in a numerically stable way during the optimization process.\

The fourth category are the \emph{frequency shift} methods.\ 
In a sense, they are a subcategory of relaxation methods applicable specifically to dynamic problems.\ 
In the context of dynamic compliance minimization, optimal designs shift the eigenfrequencies away from the loading frequency.\ 
If the loading frequency lies below an eigenfrequency, that eigenfrequency is shifted upwards by increasing the stiffness for the corresponding mode, and vice versa.\ 
This means that, for low loading frequencies, the result is a stiff design where all eigenfrequencies have been shifted upwards.\
This realization lead Olhoff and Du~\cite{olhoff2016_incremental_frequency} to propose the \emph{incremental frequency technique} for the minimization of dynamic compliance for single- and bi-material structures.\
The technique uses a continuation scheme to vary the load frequency during the optimization process and obtain better local optima.\
The generalized incremental frequency method repeats this process for a range of initial load frequencies, essentially obtaining a global optimization methodology from the first category.\
As shown by Silva et al.~\cite{Silva2019}, the underlying problem tackled in~\cite{olhoff2016_incremental_frequency} is that the dynamic compliance exhibits an alternation of resonances and anti-resonances.\ 
The resonances provide boundaries between which the optimizer gets stuck and the antiresonances provide local optima to which the optimizer converges.\ 
Guzman et al.~\cite{Guzman2023_antiresonance_matching} and Meng et al.~\cite{Meng2024_antiresonance_constraint} exploit this property by employing an antiresonance matching objective and an antiresonance constraint, respectively.\
The incremental frequency technique has also been used to gradually adapt the targeted frequency of the optimization.\ 
For example, Giannini et al.~\cite{GIANNINI2022104352} obtain MEMS designs for a target frequency of $16\, \mathrm{kHz}$ and $24\, \mathrm{kHz}$ by starting from an adapted solution found for a $26\, \mathrm{kHz}$ target.\ 
However, a rigorous analysis of this technique is lacking in the literature.\

\begin{figure}[h]
    \centering
    \includegraphics[width=\linewidth]{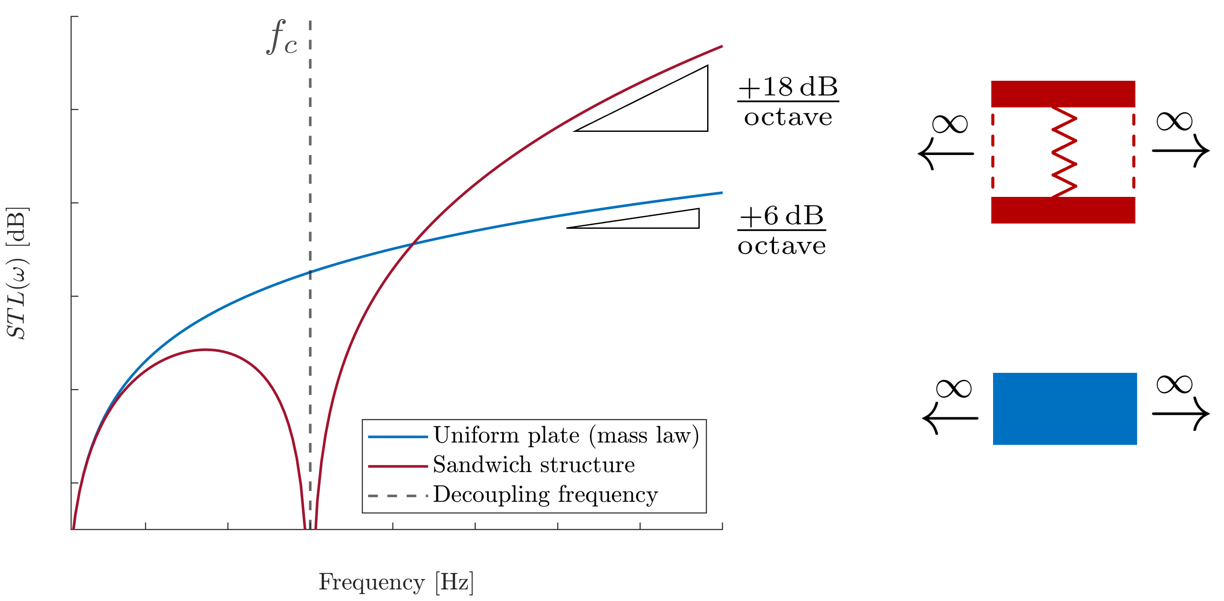}
    \caption{Sound transmission loss of a uniform plate (blue) and a sandwich panel connected with zero-mass spring (red). Both designs have the same mass.}
    \label{fig:STL_primer}
\end{figure}

The focus of this study is to reduce the number of reruns through so-called continuation strategies, i.e., optimization routines with changing hyperparameters.\ 
Such strategies are generally easy to implement and do not increase the computation time per optimization iteration.\
This work's first contribution is to develop novel continuation strategy variants for exclusion and relaxation, as well as 4 instances of the frequency shifting technique from Olhoff and Du~\cite{olhoff2016_incremental_frequency}.\ 
Global optimization is not considered due to its high computational cost.\ 
The second contribution is to compare the proposed variants through a Monte Carlo sampling of the initial guess.\ 
The third contribution consists of practical guidelines for researchers to extend and adapt the proposed continuation strategies to other settings.\\

Two applications are studied.\
The first application is the dynamic compliance minimization of the reinforced cantilever from Bendsøe and Sigmund (see \Cref{fig:Fig2.4_BendsoeSigmund}), which can yield statically stiff but dynamically low-performing designs, especially at low driving frequencies.\
Its familiarity within the topology optimization community lends this case suitable to illustrate the problem studied in this work.\
 
The second application is the maximization of a sandwich panel's sound transmission loss (STL), which measures the sound power transmitted through a structure excited by an incident sound wave.\
Although the vibro-acoustic setting is computationally demanding, it provides a useful theoretical reference known as the mass law~\cite{Norton_Karczub_2003} (see \Cref{fig:STL_primer}), which predicts that an infinite uniform plate exhibits an STL increase of $6\,\mathrm{dB}$ per octave (doubling of the frequency).\ 
In contrast, a sandwich panel, such as the 1D mass-spring-mass example of \Cref{fig:STL_primer}, achieves $18\,\mathrm{dB}/\mathrm{octave}$  above the decoupling frequency, which is a resonance of the structure.\
Any optimal solution can therefore be evaluated against the mass law because a design that does not exceed it can be replaced by a uniform plate of equal mass.\
Cool et al.~\cite{cool_TO_VA} show that STL maximization may lead to designs with mass-law performance and designs that strongly outperform it.\
For this reason, the STL maximization case is selected to quantify and compare the different continuation strategies.\\

The outline of this manuscript is as follows.\
\Cref{sec:baseline_framework_description} lays out the two cases and corresponding optimization frameworks, namely (i) the cantilever from Bendsøe and Sigmund optimized for dynamic compliance and (ii) the sandwich panel STL maximization from Cool et al.~\cite{cool_TO_VA}.\
\Cref{sec:problem_description} uses Monte Carlo sampling to illustrate and quantify the problem tackled in this work: the high probability of convergence towards low-quality optima.\
Next, \Cref{sec:strategy_investigation} elaborates on and compares the three strategy classes (exclusion, frequency shift and relaxation) by developing and implementing novel variants for each strategy class and testing them with Monte Carlo sampling.\
The first strategy (\Cref{sec:gmin_comparison}) is a novel exclusion method.\
On top of the common lower stiffness bound, the proposed continuation strategy adopts an additional upper bound on the stiffness to render low-quality optima infeasible.\
The second strategy (\Cref{sec:fcont_comparison}) is a frequency shift method, i.e., a variant of the incremental frequency technique from Olhoff and Du~\cite{olhoff2016_incremental_frequency}, that pulls high-frequency solutions to lower frequency ranges.\
The study shows that, although very beneficial, frequency continuation suffers from high computational cost and imperfect success rates when applied to the sandwich panel case.\ 
The third strategy (\Cref{sec:robustness_comparison}) is a relaxation method, whose main novelty is a softening of the design robustness at the start of the optimization process.\
\Cref{subsec:comparison_strategies} bundles the conclusions for each strategy, \Cref{subsec:generality_and_guidelines} discusses their applicability and extension to other contexts.\
Finally, \Cref{sec:conclusions} formulates the conclusions.

\section{Optimization framework description} \label{sec:baseline_framework_description}

This section lays out the two cases. Focus lies on boundary conditions and optimization setup, with numerical details, material and geometric properties in \ref{app:TO}.

\subsection{Optimization variables}

Both the cantilever and the sandwich panel core are optimized with density-based topology optimization.\
That is, each element is controlled with an optimization variable $\bm{\xi}$ and projection and filtering techniques map the optimization variable field $\bm{\xi}$ to the physical density field $\bm{\xi}_P$ used to evaluate material properties.\ 
For the cantilever, $\xi_{i,P} \in [0,1]$ interpolates between void ($\xi_{i,P}=0$) and solid material  ($\xi_{i,P}=1$).\
Similarly, for the sandwich panel core, $\xi_{i,P} \in [0,1]$ interpolates between fluid ($\xi_{i,P}=0$) and solid material  ($\xi_{i,P}=1$).\
The design representation for both cases uses the double filtering approach of Christiansen et al.~\cite{christiansen2015doublefilt}, with filter radii $R_1$ and $R_2$ and corresponding smoothed Heaviside projections with steepness variables $\beta_1$ and $\beta_2$.\
Furthermore, the robust formulation of~\cite{wang2011projection} is used, yielding three physical designs: an eroded ($e$), a blueprint ($b$) and a dilated ($d$) design.\ 
The minmax formulation, discussed further below, ensures the optimizer prioritizes the worst-performing of the three designs.\\ 

The cantilever uses the augmented PDE filter of Wallin et al.~\cite{Wallin2020} to enforce the length scale and prevent the design from sticking to the boundaries.\ 
The sandwich panel case uses a periodic version of the convolutional filter with a hat function kernel~\cite{bourdin2001filters} to enforce length scale and periodicity.\

\subsection{Dynamic compliance minimization of the cantilever from Bendsøe and Sigmund}

\Cref{fig:setup_BS} shows the design domain and boundary conditions of the cantilever, which has its outer frame fixed at solid material.\ 
The structure is fixed at the left and loaded with a vertical, time-periodic force $F\cos(\omega t)$, with $\omega=2\pi f$ and $f$ the driving frequency.\
Finite element discretization yields a linear system in the Fourier domain

\begin{equation}
\left( (1+\alpha \mathrm{i})\mathbf{K} -\omega^2 \mathbf{M} \right) \mathbf{u} = \mathbf{f},
\end{equation}

where $\mathbf{K}$ and $\mathbf{M}$ are the stiffness and mass matrices, $\mathbf{u} $ and $\mathbf{f}$ are the displacement and force vectors and $\alpha$ is the structural damping coefficient.\ 
The objective function is the logarithm of the dynamic compliance $c_{b/e/d} = \log_{10}(|\mathbf{u}^T \mathbf{f}|)$, where $|\cdot|$ expresses the complex magnitude and the $b/e/d$ underscript denotes whether the blueprint, eroded or dilated design is considered.\
For numerical stability, the objective function is normalized:
\begin{linenomath}
\begin{equation}
\label{eq:dynamic_compliance_normalization}
    J_{b/e/d} = \frac{c_{b/e/d}(\omega)}{A}+B,
\end{equation}
\end{linenomath}
where $A=|c_{\mathrm{c_{\text{solid}}}}|$ with $c_{\text{solid}}$ the dynamic compliance of the fully solid design and $B=2$ such that the objective function starts at a value close to $1$.\ 
Since no normalization issues were encountered during this study, the normalization was not tuned further.

\subsection{STL maximization of a sandwich panel}

\begin{figure}[!t]
\centering
\begin{subfigure}[t]{0.55\textwidth}
  \centering
\includegraphics[width=\textwidth]{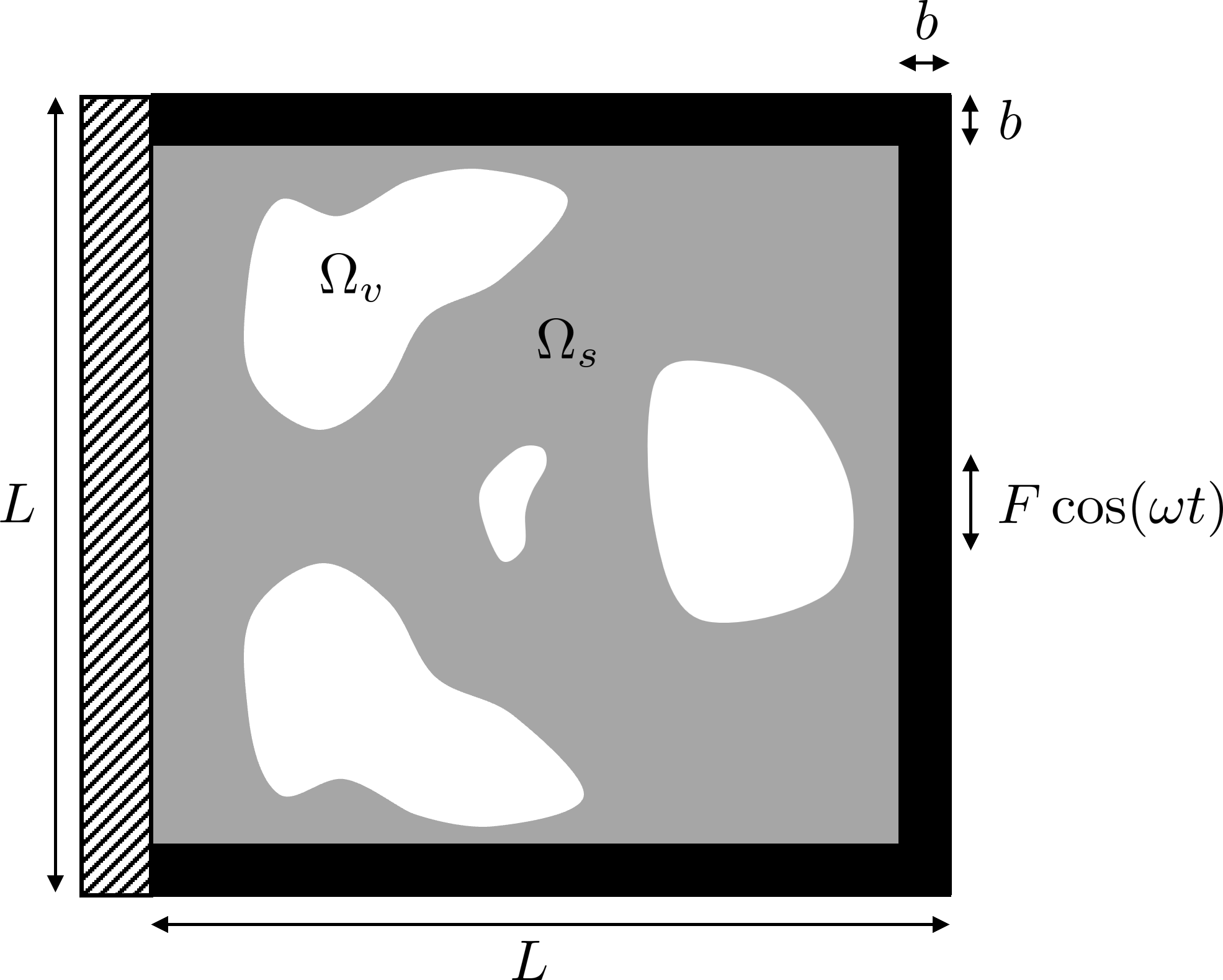}
  \caption{Reinforced cantilever.} \label{fig:setup_BS}
\end{subfigure}%
\begin{subfigure}[t]{0.45\textwidth}
  \centering
\includegraphics[width=\textwidth]{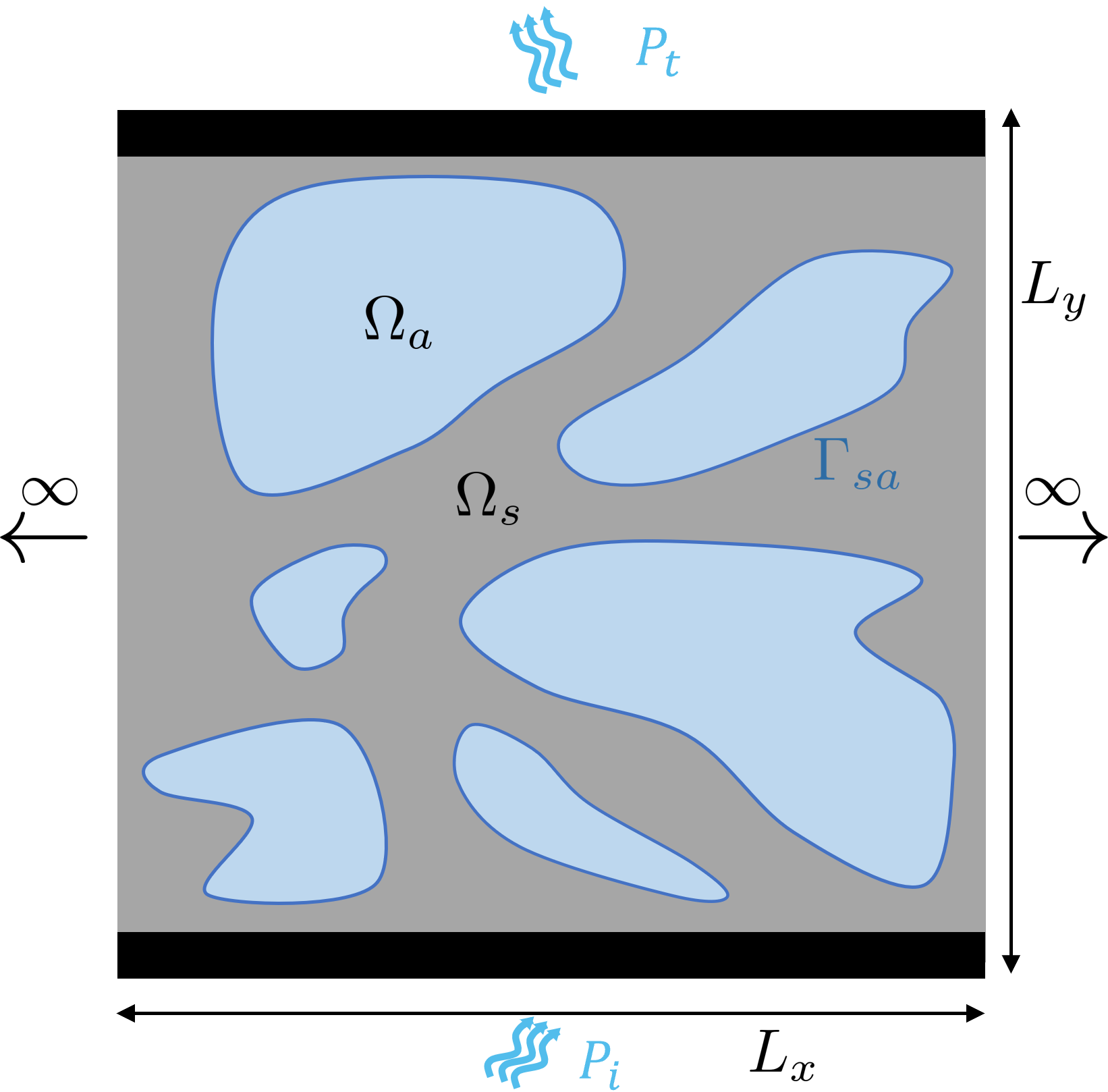}
  \caption{Sandwich panel.} \label{fig:scheme_prob}
\end{subfigure}
\caption{Boundary conditions and design domain of the considered cases.}
\end{figure}

Fig.~\ref{fig:scheme_prob} shows the sandwich panel setup. It consists of two plates with a vibro-acoustic core in between.\
The core contains an infinitely repeating unit cell. \

Above and below the panel are infinite acoustic half spaces containing air (density $\rho_a=1.225$~kg/m$^3$ and speed of sound $c_a=340$~m/s).\
An incident plane acoustic pressure wave (amplitude $P_i$, acoustic wave number $k_a=\omega/c_a$ and radial frequency $\omega$ (rad/s)) excites the panel at the bottom.\ 
This causes a transmitted wave at the top (amplitude $P_t$) and a reflected wave at the bottom (amplitude $P_r$).\
In the Fourier domain, these three plane acoustic pressure waves are written as:
\begin{linenomath}
\begin{equation}
\label{eq:Pi_Pt_Pr}
\begin{aligned}
    p_i(x,y) = P_i \mathrm{e}^{\mathrm{i}(-k_xx-k_yy)}, \\
    p_t(x,y) = P_t \mathrm{e}^{\mathrm{i}(-k_xx-k_yy)}, \\
    p_r(x,y) = P_r \mathrm{e}^{\mathrm{i}(-k_xx+k_yy)},
\end{aligned}
\end{equation}
\end{linenomath}
where $k_x=k_a sin(\theta)$, $k_y=k_a cos(\theta)$ are trace wavenumbers and $\theta$ is the incidence angle. For ease of notation, the temporal factor $\mathrm{e}^{\mathrm{i}\omega t}$ is excluded.\
The STL is used to asses the acoustic performance of the structure. It is defined as the ratio of the transmitted pressure amplitude to the incident pressure amplitude~\cite{fahy2007sound}: 
\begin{linenomath}
\begin{equation}
    \label{eq:STL}
\begin{aligned}
    \mathrm{STL(\omega,\theta)} &= -10\log_{10}{\left(\tau(\omega,\theta) \right)} \\
    \tau(\omega,\theta) &= \left| \frac{P_t}{P_i}\right|^2.
\end{aligned}
\end{equation}
\end{linenomath}
Here, $\tau$ denotes the sound power transmission coefficient.\
Finite Element discretization yields a linear system (see \ref{app:TO}) that must be solved for the nodal displacement and pressure.\
Bloch-Floquet boundary conditions~\cite{bloch1929quantenmechanik} are imposed on the left and right edges of the unit cell and a force on the bottom edge to  model the acoustic half-spaces~\cite{yang2019vibroacoustic}.\
For a circular frequency $\omega$, the STL is computed with the WFE method~\cite{yang2017prediction,cool_TO_VA,boukadia2020wave}.\
A range $\Delta \omega = [\omega_{-}, \omega_{+}]$ is considered by numerically approximating the average STL as follows:
\begin{linenomath}
\begin{equation}
\label{eq:STL_averaging}
\begin{aligned}
    STL(\Delta\omega, \theta)&=\frac{\int^{\omega_{+}}_{\omega_{-}}STL(\omega, \theta) \,\mathrm{d}\omega}{\int^{\omega_{+}}_{\omega_{-}}\,\mathrm{d}\omega}\\ &\approx \frac{1}{\omega_{+}-\omega_{-}}\sum_{\omega_i \in \Delta\omega} a_{\omega_i}STL(\omega_i, \theta),
\end{aligned}
\end{equation}
\end{linenomath}
where $a_{\omega_i}$ are (in this work) five midpoint integration weights.\ 
The incidence angle $\theta$ is kept at $0\degree$ and, for brevity, omitted in what follows.\ 
In the remainder of this manuscript, all numerical values of $\omega$ are specified in $\mathrm{Hz}$ without the implied $2\pi$ scaling factor.\
Finally, for numerical stability, the STL is normalized:
\begin{linenomath}
\begin{equation}
\label{eq:STL_normalization}
    J_{b/e/d} = \frac{STL_{b/e/d}(\Delta\omega)}{C},
\end{equation}
\end{linenomath}
where $b/e/d$ denotes either the blueprint, eroded or dilated design and $C=120$ is a problem-dependent normalization constant to ensure the objective starts close to $1$ (see Cool et al.~\cite{cool_TO_VA}).\

\subsubsection{Constraints} \label{sec:constraints}

Both cases employ three constraints: 1) a bound constraint on the optimization variables $\bm{\xi}$, 2) a volume constraint and 3) a connectivity constraint.\
The bound constraint restricts the variables to lie in $[0, 1]$, so $0 \leq \bm{\xi} \leq 1$.\
The volume constraint imposes a maximum volume of $V=0.5V_{tot}$, with $V_{tot}$ the total volume of the design domain, and is applied on the dilated design.\
It is added for similar reasons as in the work by Bendsøe and Sigmund~\cite{BendsøeSigmundBook2003}: without it, routines often converge to fully solid designs.
The connectivity constraint consists of a self-weight load to reduce the appearance of disconnected structures in the core~\cite{cool_connectivity_review}.\
The compliance under this self-weight $\theta_{sw}$ is constrained to be $\mu_{sw}$ times smaller than the self-weight compliance for a fully solid design, denoted by $\hat{\theta}_{sw}$.\
For both cases, the connectivity constraint is applied on the eroded design.\\
As will be shown later, despite the connectivity constraint, a minority of the designs in this study exhibit disconnected regions.\
More specifically, some beams contain greyscale density values since the optimizer aims to obtain a specific stiffness to mass ratio and hence control its resonance frequency.\
For the cantilever, this is resolved by ending the optimization with an explicit penalization of the greyness (see \ref{app:TO}).\
This is not done for the sandwich panel core to reduce the number of iterations.\ 
The effect on the conclusions of this study is limited: if perfect connectivity is enforced, the $\mathrm{STL}$ of some high-performing sandwich panels is slightly reduced but the underlying mechanisms still outperform the mass law by a considerable margin.\
As perfect connectivity is not the focus, this is not further explored.\

\subsubsection{Formulation}

The baseline implementation uses the following minmax formulation:
\begin{linenomath}
\begin{equation}
\label{eq:minmax_formulation}
\left\{
\begin{aligned}
\min_{\bm{\xi}\in \mathbb{R}^{N_e}} \quad \max (&J_b, J_e, J_d) \\
\textrm{s.t.} \quad \quad \quad \quad \quad J_{\mathrm{vol}}=~& v_{d,\mathrm{P}}/V-1  \leq 0\\
 J_{\mathrm{conn}} =~& \theta_{\mathrm{sw}}/\mu_{\mathrm{sw}}\hat{\theta}_{\mathrm{sw}} - 1 \leq 0\\
 0 \leq \bm{\xi} &\leq 1.\\
\end{aligned}
\right.
\end{equation}
\end{linenomath}

Here, the volume and connectivity constraints are normalized to obtain $J_{\mathrm{vol}}$ and $J_{\mathrm{conn}}$, respectively.\ 
The minmax formulation is solved with the Method of Moving Asymptotes~\cite{svanberg1987method} (MMA), where the slack variable $z$ circumvents the non-differentiable $\max$ operator.  

\section{Problem description} \label{sec:problem_description}

This section illustrates (\Cref{sec:problem_illustration}) and quantifies (\Cref{sec:problem_quantification}) the problem tackled in this work.\

\subsection{Problem illustration} \label{sec:problem_illustration}

\begin{figure}[h]
\centering
\includegraphics[width=0.5\textwidth]{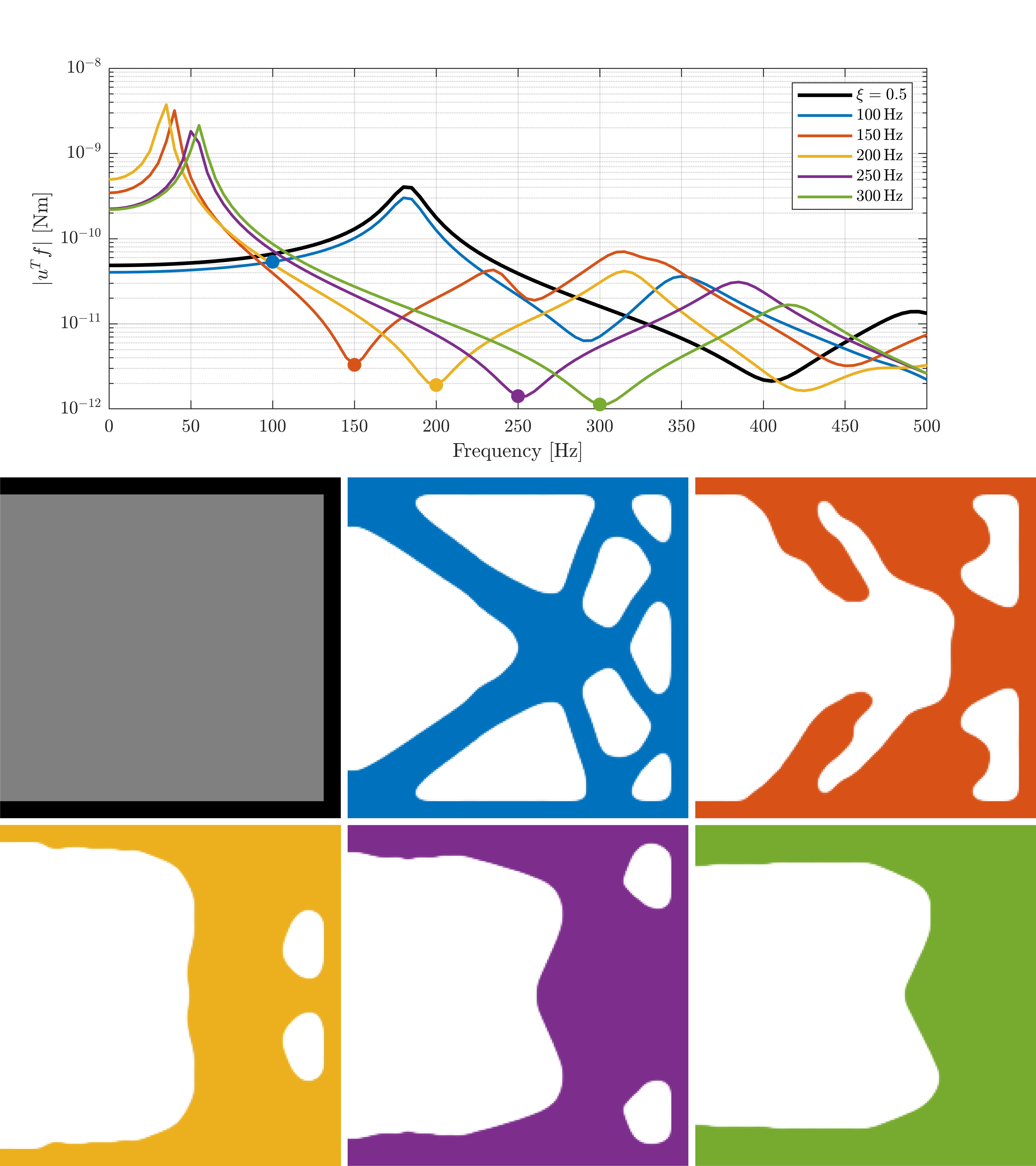}
\caption{Solutions of the cantilever problem for various driving frequencies. Top: dynamic compliance as function of the frequency. Bottom: uniform initial guess and corresponding blueprint designs.}
\label{fig:BendsoeSigmundComparison}
\end{figure}

\Cref{fig:BendsoeSigmundComparison} shows solutions of the cantilever case using the framework described in the previous section, for different driving frequencies.\
The first natural frequency of the uniform initial guess lies around $f_1=180$~Hz.\
The optimization run with driving frequency $100$~Hz, denoted in blue, results in a statically stiff solution much like \Cref{fig:Fig2.4_BendsoeSigmund}a-b.\
In contrast, for only slightly higher driving frequencies, high-performing designs are obtained.\ 
This shows a large dependence on the initial guess and the driving frequency.\

\begin{figure}[!htb]
\centering
\begin{subfigure}[t]{0.5\textwidth}
  \centering
\captionsetup{width=0.95\textwidth}
  \includegraphics[width=0.95\textwidth]{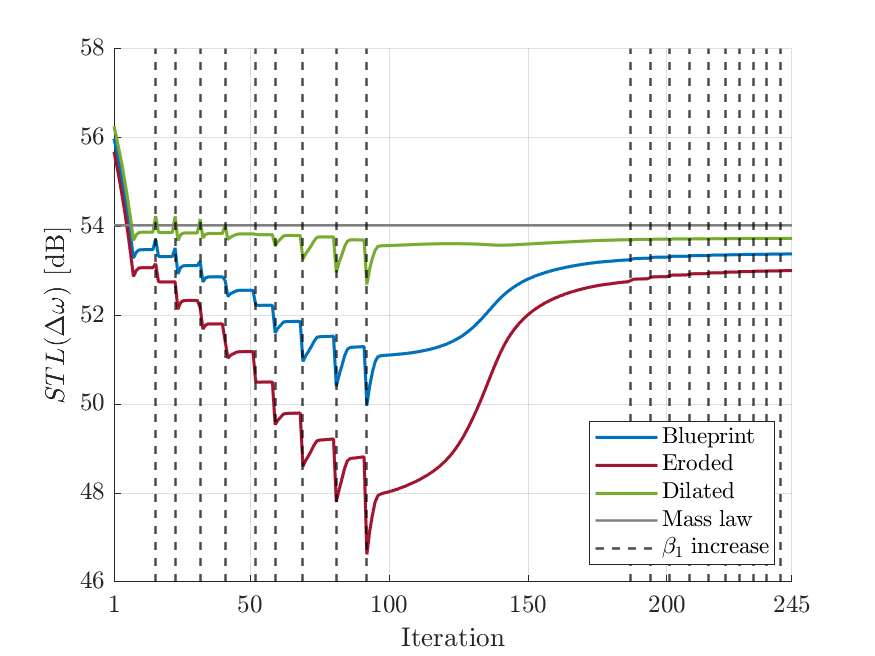}
  \caption{Sound transmission loss $\mathrm{STL}(\Delta \omega)$.} \label{subfig:STL_static}
\end{subfigure}%
\begin{subfigure}[t]{0.5\textwidth}
  \centering
  \captionsetup{width=0.95\textwidth}
  \includegraphics[width=0.95\textwidth]{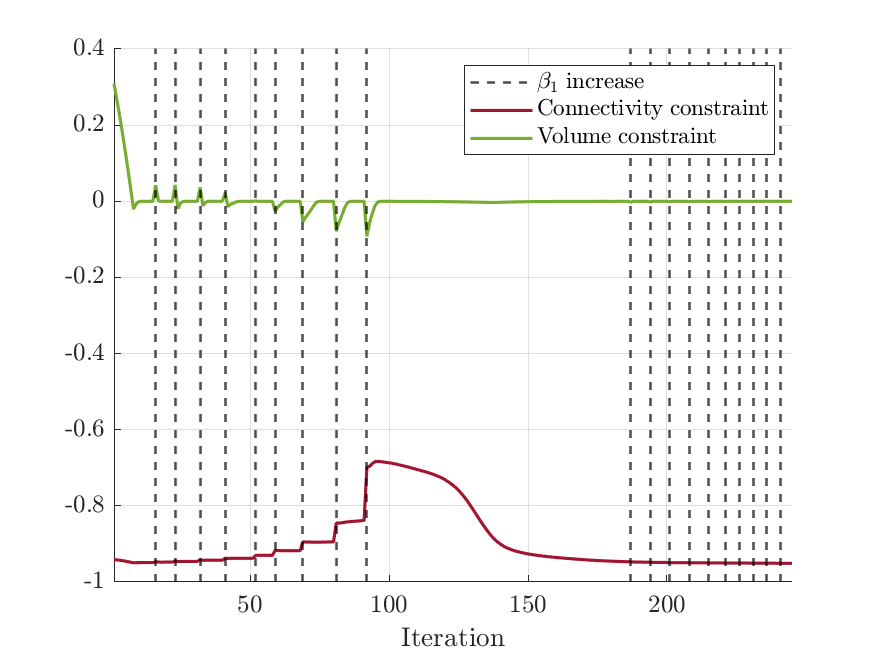}
  \caption{Connectivity and volume constraints.} \label{subfig:cg_static}
\end{subfigure}\\
\begin{subfigure}[c]{0.5\textwidth}
  \centering
    \captionsetup{width=0.95\textwidth}
  \includegraphics[width=0.95\textwidth]{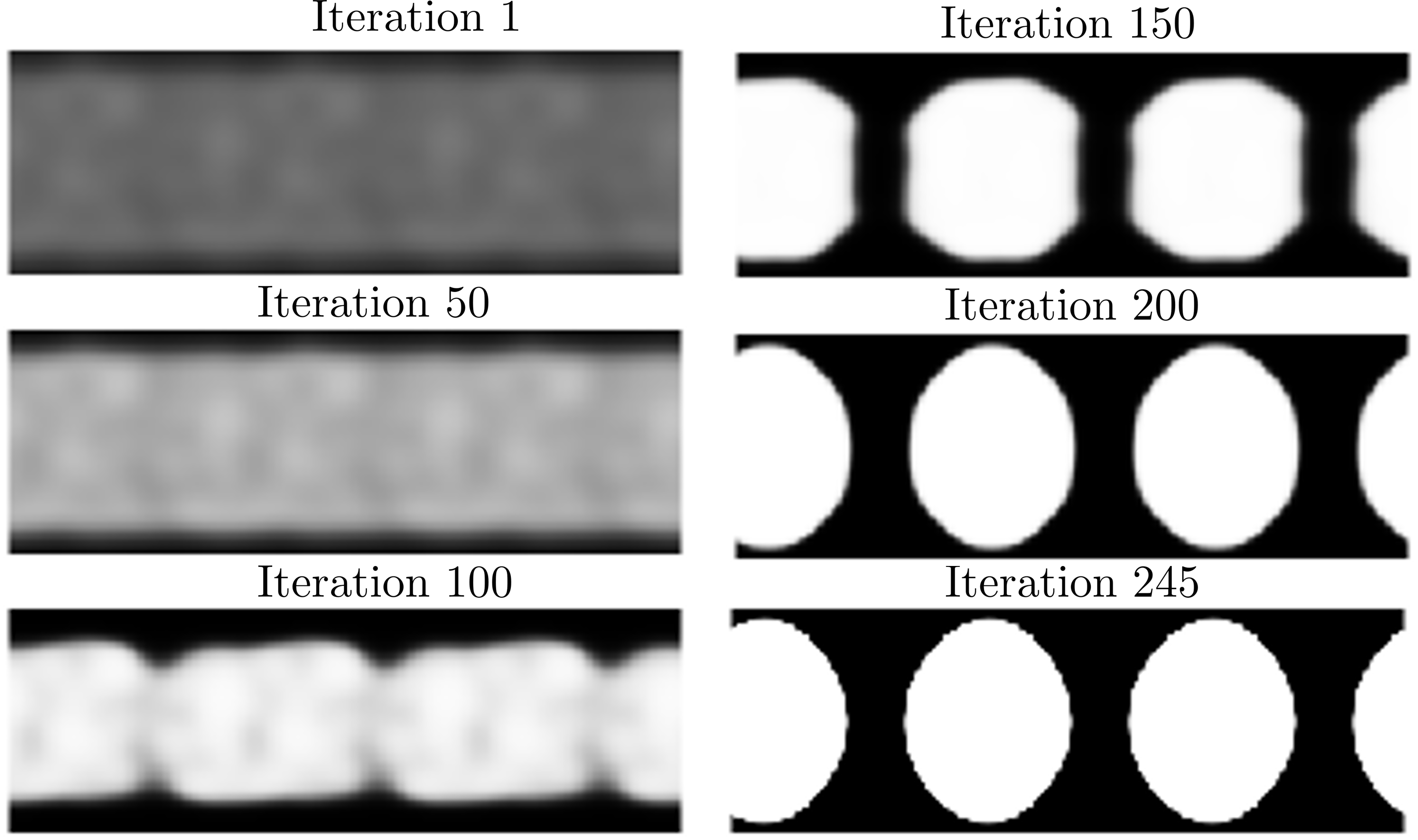}
  \caption{Blueprint designs, repeated 3 times.} \label{subfig:designs_static_run}
\end{subfigure}%
\begin{subfigure}[c]{0.5\textwidth}
  \centering
  \captionsetup{width=0.95\textwidth}
  \includegraphics[width=0.8\textwidth]{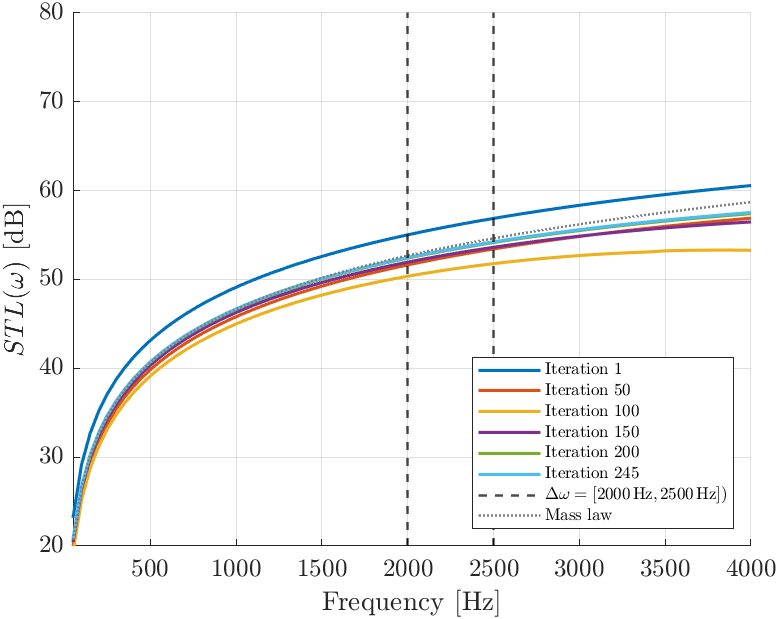}
  \caption{Sound transmission loss $\mathrm{STL}(\omega)$ of blueprint designs.}
  \label{subfig:fullSTL_static}
  \end{subfigure}%
\caption{Baseline optimization run yielding a mass-driven optimum, for $\Delta \omega=[2000\,\mathrm{Hz},2500\,\mathrm{Hz}]$.}
\label{fig:problem_illustration_static}
\end{figure}

\begin{figure}[!htb]
\centering
\begin{subfigure}[t]{0.5\textwidth}
  \centering
\captionsetup{width=0.95\textwidth}
  \includegraphics[width=0.95\textwidth]{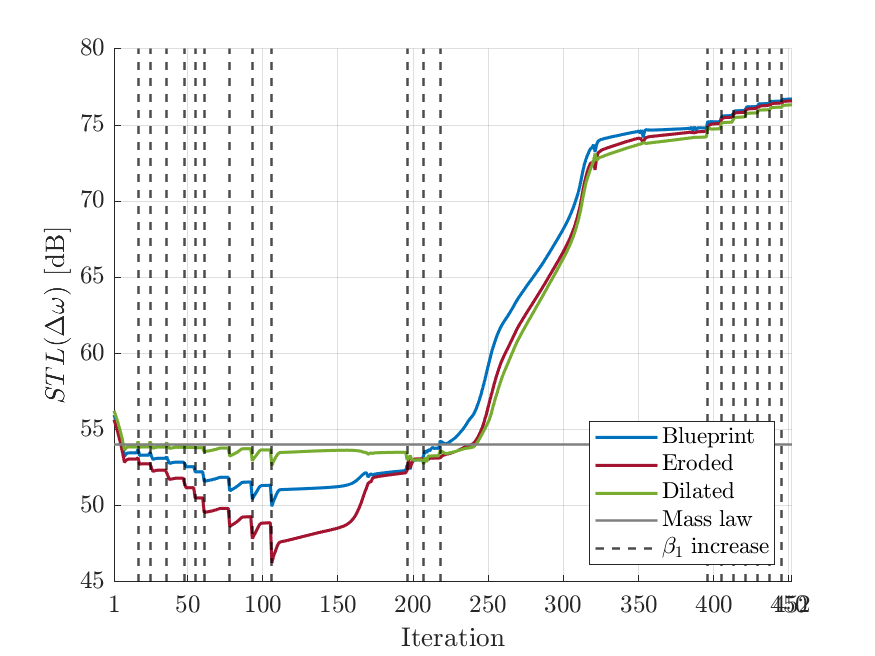}
  \caption{Sound transmission loss $\mathrm{STL}(\Delta \omega)$.}  \label{subfig:STL_dynamic}
\end{subfigure}%
\begin{subfigure}[t]{0.5\textwidth}
  \centering
  \captionsetup{width=0.95\textwidth}
  \includegraphics[width=0.95\textwidth]{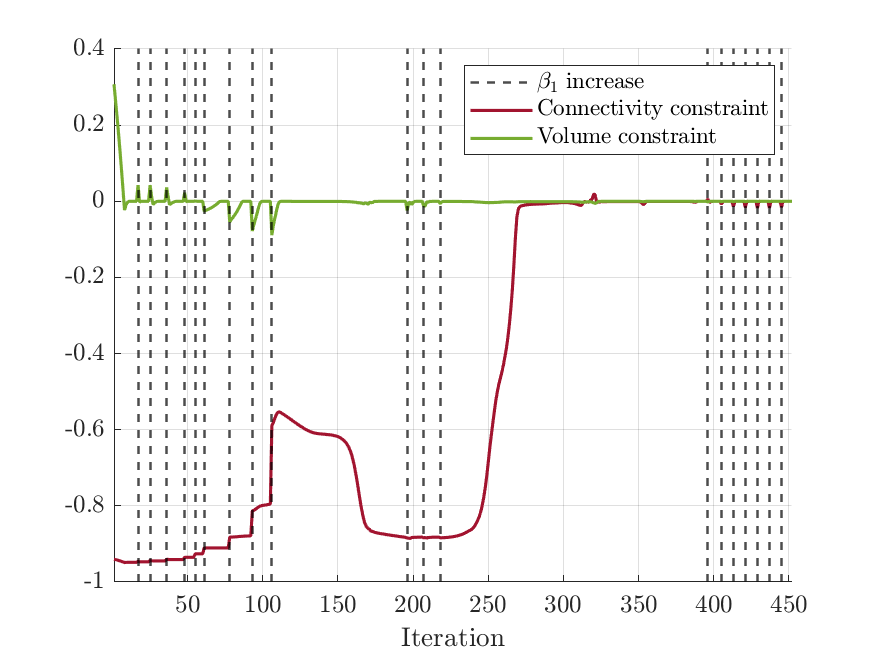}
  \caption{Connectivity and volume constraints.}\label{subfig:cg_dynamic}
\end{subfigure}\\
\begin{subfigure}[c]{0.5\textwidth}
  \centering
    \captionsetup{width=0.95\textwidth}
  \includegraphics[width=0.95\textwidth]{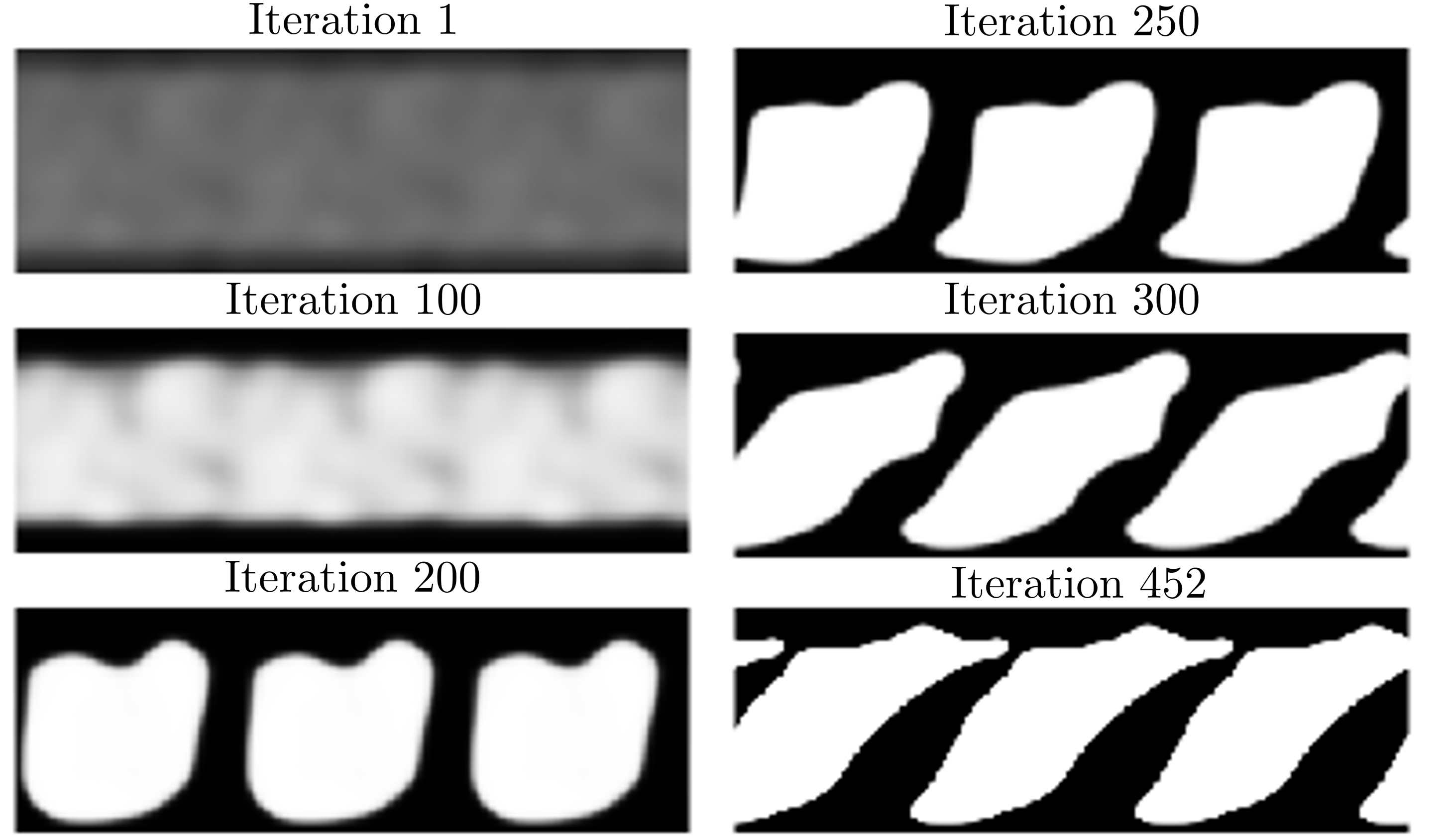}
  \caption{Blueprint designs, repeated 3 times.}\label{subfig:designs_dynamic_run}
\end{subfigure}%
\begin{subfigure}[c]{0.5\textwidth}
  \centering
  \captionsetup{width=0.95\textwidth}
  \includegraphics[width=0.8\textwidth]{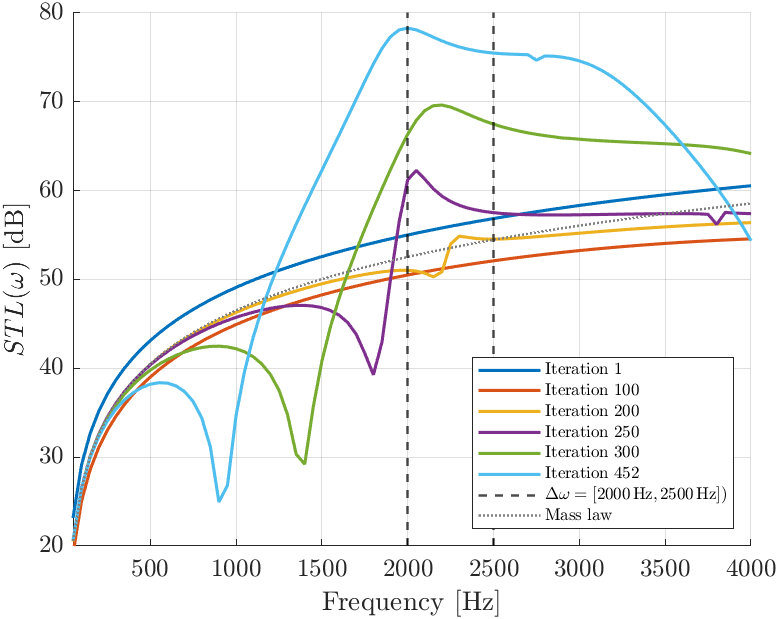}
  \caption{Sound transmission loss $\mathrm{STL}(\omega)$ of blueprint designs.}
  \label{subfig:fullSTL_dynamic}
  \end{subfigure}%
\caption{Baseline optimization run yielding a high-performing optimum, for $\Delta \omega=[2000\,\mathrm{Hz},2500\,\mathrm{Hz}]$.}
\label{fig:problem_illustration_dynamic}
\end{figure}

\Cref{fig:problem_illustration_dynamic,fig:problem_illustration_static} show two optimization runs for the sandwich panel optimization, targeting the frequency range $[2000\,\mathrm{Hz},2500\,\mathrm{Hz}]$ by starting from a random initial guess.\
A clear difference is visible.\
Whereas \Cref{fig:problem_illustration_dynamic} shows a successful run that achieves a high-performing design, \Cref{fig:problem_illustration_static} shows a run that remains stuck in a low-performing local optimum that behaves like a single plate, i.e., with mass law peformance.\

The unsuccessful run starts from a random initial guess (\Cref{subfig:designs_static_run}, top left). The first $100$ iterations follow a repeating pattern: (i) a violation of the volume constraint (due to a $\beta_1$ increase, see \Cref{subfig:cg_static}), (ii) recovery to a feasible design, (iii) an $\mathrm{STL}(\Delta \omega)$ (\Cref{subfig:STL_static}) that flatlines and finally (iv) a $\beta_1$ increase that again violates the volume constraint.\ 
During this process, the core of the sandwich panel loses stiffness, shifting the decoupling frequency down and hence causing an $\mathrm{STL}(\Delta \omega)$ decrease.\ 
In the remaining $145$ iterations, the optimize counters this effect: the panel is connected with a rigid, vertical beam that increases the core stiffness, shifts the decoupling frequency upwards and increases the $\mathrm{STL}(\Delta \omega)$.\ 
However, the final performance is unsatisfactory.\

The successful run starts with a different random initial guess.\ 
The first $150$ iterations see the same softening of the core followed by a stiffening due to the formation of a rigid, vertical connection.\ 
After iteration $150$, the run visibly diverges: the vertical connection becomes askew (\Cref{subfig:designs_dynamic_run}) until the connectivity constraint is active (\Cref{subfig:cg_dynamic}).\ 
Furthermore, the decoupling frequency drops $> 1000\,\mathrm{Hz}$ below the targeted frequency range (\Cref{subfig:fullSTL_dynamic}) and the performance ($\mathrm{STL}(\Delta \omega)$, \Cref{subfig:STL_dynamic}) increases by $20\,\mathrm{dB}$.\ 
Around iteration $400$, the increase flatlines and subsequent $\beta_1$ increases result in a crisp, high-performing design.\
The high performance is due to two mechanisms: 1) a reduction in the decoupling frequency to exploit the strong $18\,\mathrm{dB}$ per octave STL increase and 2) the transformation of the vertical motion of the top plate, which causes the transmitted pressure wave, to a horizontal-dominated motion that does not substantially excite the acoustic domain.\
As a result, the STL of the high-performing design lies $20\,\mathrm{dB}$ or about $35\%$ above the mass law in the targeted frequency range.\ 
This is much better than the low-performing, mass-driven optimum of \Cref{fig:problem_illustration_static}.\ 

\subsection{Problem quantification} \label{sec:problem_quantification}

The previous section illustrates a large dependence on the initial guess and frequency: minor deviations result in large performance differences.
To quantify this behavior, this section employs a Monte Carlo sampling methodology to compute the probability of achieving a high-performing optimum ($P_{\text{HP}}$).\\

 For the cantilever, the optimization routine was run for 11 frequencies, from $50\,\mathrm{Hz}$ to $300\,\mathrm{Hz}$, with $20$ random initial guesses.\
The driving frequencies were chosen to lie both above and below the first natural frequency.\
\ref{app:monte_carlo_sampling_methodology} elaborates upon the choice of initial guesses.\

\begin{figure}[h]
\centering
\includegraphics[width=\linewidth]{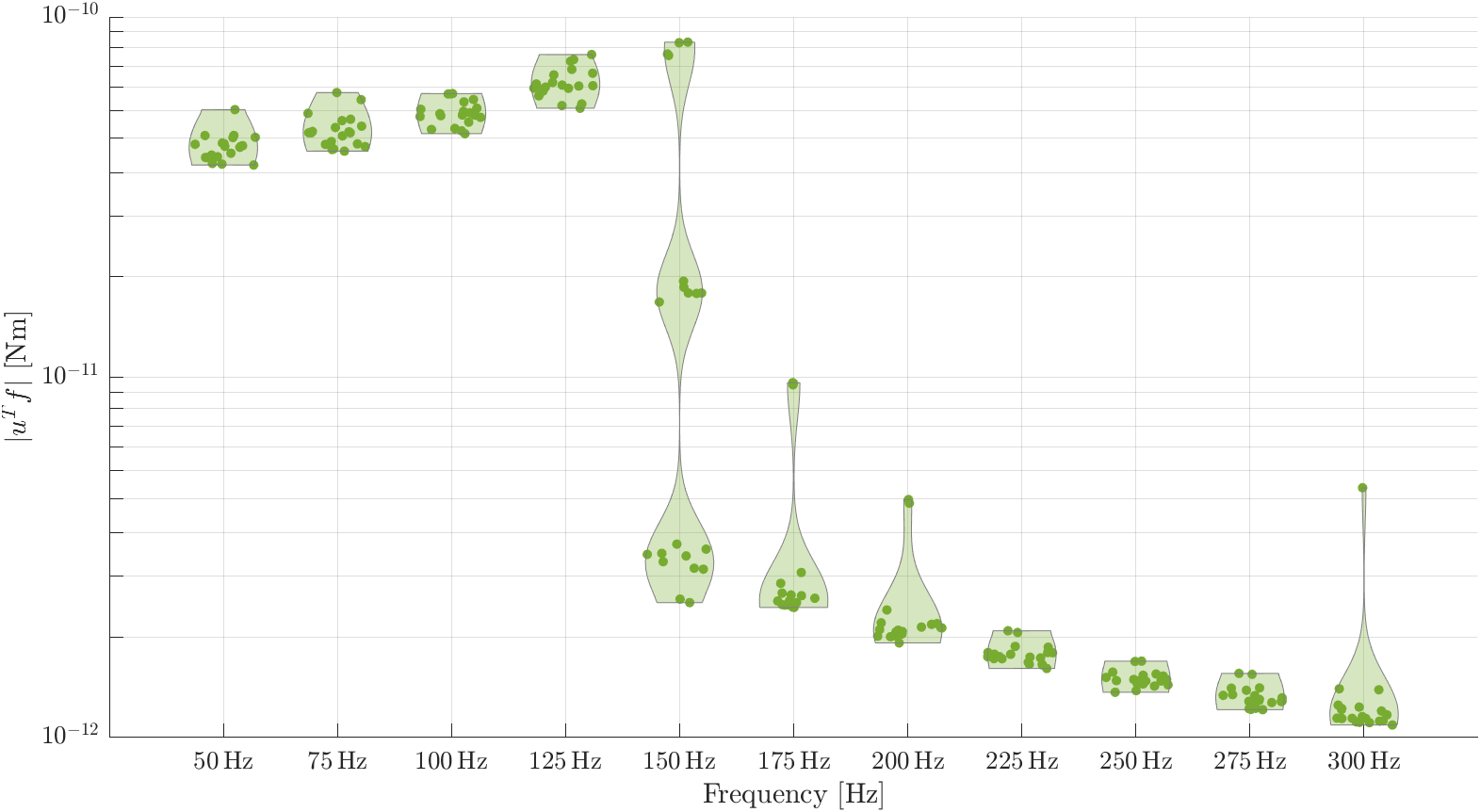}
\caption{Monte Carlo sampling methodology applied to the cantilever. The driving frequency ($x$-axis) is plotted against the performance $\log(|\mathbf{u}^T\mathbf{f}|)$  of the blueprint design ($y$-axis). Violin plots illustrate sample density and contain small, artificial frequency shifts to enhance visibility~\cite{violin}.}
\label{fig:BendsoeSigmund_baseline}
\end{figure}

\Cref{fig:BendsoeSigmund_baseline} plots the obtained performances as a function of the frequency for the cantilever.\ 
For frequencies $125\,\mathrm{Hz}$ and below, only statically stiff designs were obtained.\
For  frequencies $175\,\mathrm{Hz}$ and above, all designs exploit an antiresonance and have good performance.\
However, $150\,\mathrm{Hz}$ is a transition point: both high- and low-performing designs are found.\
It must be stressed that designs with a good performance exist at and below $150\,\mathrm{Hz}$: Olhoff and Du~\cite{olhoff2016_incremental_frequency} define high and low frequencies as those above and below the first natural frequency and use the incremental frequency technique to draw high-frequency solutions to low frequencies.\

A similar study is now performed for the sandwich panel, revealing an additional problem at high frequencies.\
For the sandwich panel, $16$ frequency ranges $\Delta \omega$ between $50\,\mathrm{Hz}$ and $8000\,\mathrm{Hz}$ are considered, each with a width of $500\mathrm{Hz}$ except the first, $[50\mathrm{Hz},500\mathrm{Hz}]$, to exclude $0\,\mathrm{Hz}$. For each frequency range, the optimization is run $N=20$ times from a random initial guess $\bm{\xi}_i$, $i \in [1, \ldots, 20]$ (see \ref{app:monte_carlo_sampling_methodology}).\ 
Afterwards, a design's performance $\mathrm{STL}^*=\min(STL_b(\Delta \omega), STL_e(\Delta \omega), STL_d(\Delta \omega))$ is compared to the mass law $\mathrm{STL}_{ml}(\Delta \omega)$.\ 
If $\mathrm{STL}^* \geq 1.1 \times STL_{ml}$, it is considered to be \emph{high-performing}, and if $\mathrm{STL}^* < 1.1 \times STL_{ml}$, it is \emph{low-performing}.\ 
Note that although the cutoff of $10\%$ above the mass law is arbitrary, results show that the classification is relatively insensitive to the  cutoff value.\\

For a frequency range $\Delta \omega$, the number of high-performing optima, $N_{\text{HP}}$, is compared to the number of low-performing optima, $N_{\text{LP}}$, to provide an estimate $\tilde{P}_{\text{HP}}$ of the probability of achieving a high-performing optimum:

\begin{linenomath}
\begin{equation}
\label{eq:chance_of_high_performance}
P_{\text{HP}} \approx \tilde{P}_{\text{HP}} =  100\frac{N_{\text{HP}}}{N_{\text{LP}}+N_{\text{HP}}}\,\mathrm{[\%]}.
\end{equation}
\end{linenomath}

Since only $20$ samples are used, $\tilde{P}_{\text{HP}}$ has relatively high error bars.\
It is however sufficiently accurate to reveal trends.\ 
To this end, frequency ranges $\Delta \omega$ are classified as follows:
\begin{itemize}
    \item $\tilde{P}_{\text{HP}} \leq ~20\%$: low-performing region,
    \item $\tilde{P}_{\text{HP}} \geq ~80\%$: high-performing region,
    \item $\tilde{P}_{\text{HP}} \in ]~20\%,~80\%[$: transition region.
\end{itemize}

\begin{figure}[h]
\centering
\captionsetup{width=0.95\textwidth}
\includegraphics[width=0.95\textwidth]{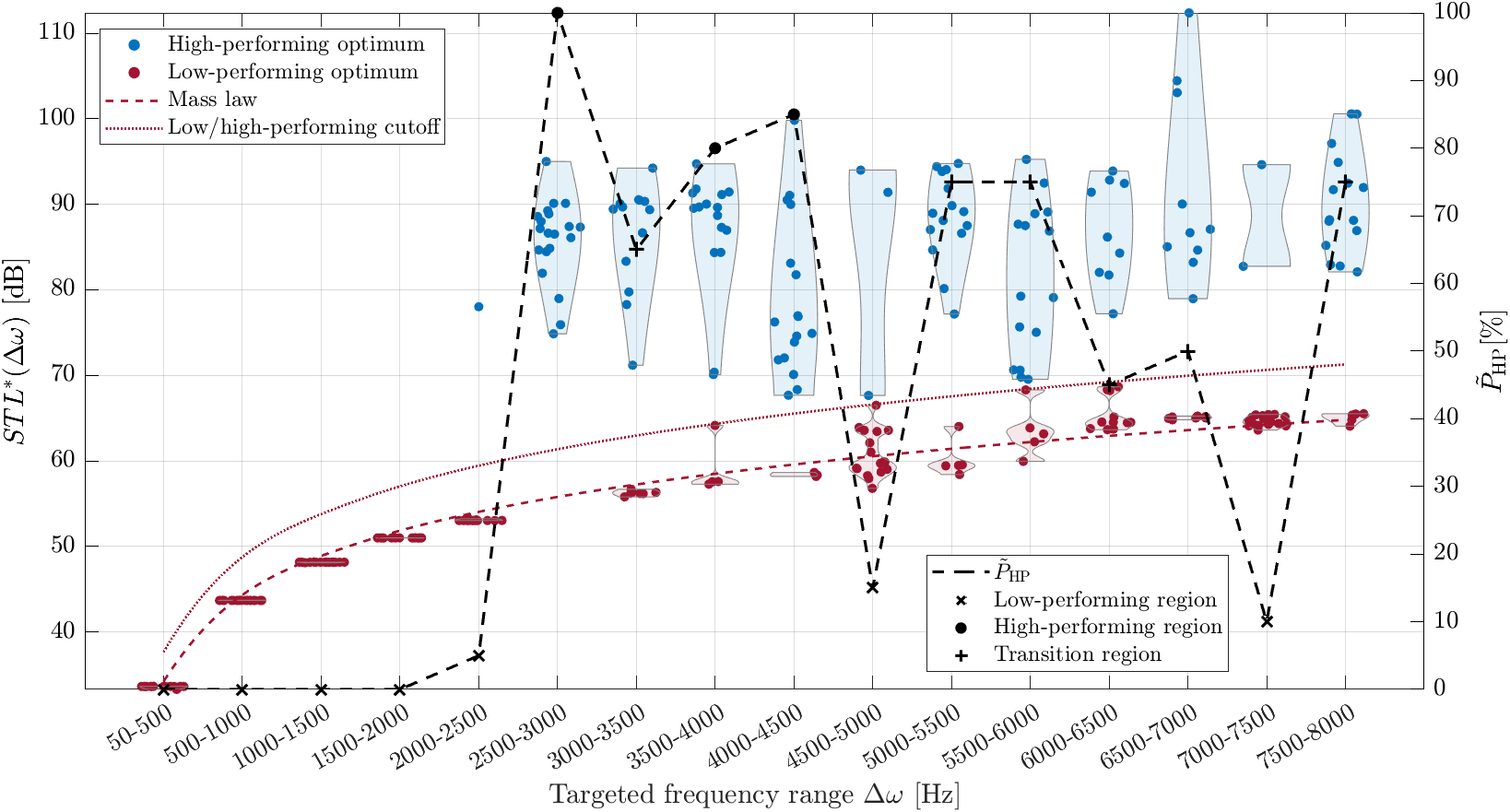}
\caption{Monte Carlo sampling methodology applied to sandwich panel case. Performances $\mathrm{STL}^*$ of the optimized designs (left $y$-axis) and estimated chance at a high-performing optimum (right $y$-axis) as a function of the frequency range $\Delta\omega$ ($x$-axis).}
\label{fig:results_baseline}
\end{figure}

\Cref{fig:results_baseline} shows the sandwich panel performance and the estimated probability $\tilde{P}_{\text{HP}}$ as a function of the frequency range.\
Two problems appear.\
First, there exists a low-performing zone at low frequencies ($\Delta\omega < 2500\,\mathrm{Hz}$), where \emph{all} frequency ranges are low-performing.\
This zone is similar to the cantilever case, where below $125$~Hz only statically stiff optima occur.\
For the sandwich panel case, the width of this zone is likely related to the decoupling frequency.\
Second, for $\Delta\omega > 4500\,\mathrm{Hz}$, $\tilde{P}_{\text{HP}}$ frequently dips below $50\%$, which contrasts strarkly with the cantilever case, which has $\tilde{P}_{\text{HP}}=100\%$ above $175$~Hz.\
Three frequency ranges are especially notable for the sandwich panel case: right after the low-performing zone, $\Delta\omega=[2500\,\mathrm{Hz}, 3000\,\mathrm{Hz}]$ attains $\tilde{P}_{\text{HP}}=100\%$, while $\Delta\omega=[4500\,\mathrm{Hz}, 5000\,\mathrm{Hz}]$ and $\Delta\omega=[7000\,\mathrm{Hz}, 7500\,\mathrm{Hz}]$ are isolated low-performing regions.\ 
This signifies the existence of, numerically speaking, highly attractive families of optima around these frequencies.\
Note that the definition of ``low' and ``high'' frequency for the sandwich panel is different from the definition of Olhoff and Du for the cantilever.\
Low frequency for the sandwich panel means $\Delta\omega < 2500\,\mathrm{Hz}$ since this empirically resembles the low-frequency region for the cantilever.\
The first natural frequency of the initial sandwich panel designs lie more than twice as high, however.\

\begin{figure}[h!]
\centering
\begin{subfigure}[t]{\linewidth}
  \centering
  \includegraphics[width=\linewidth]{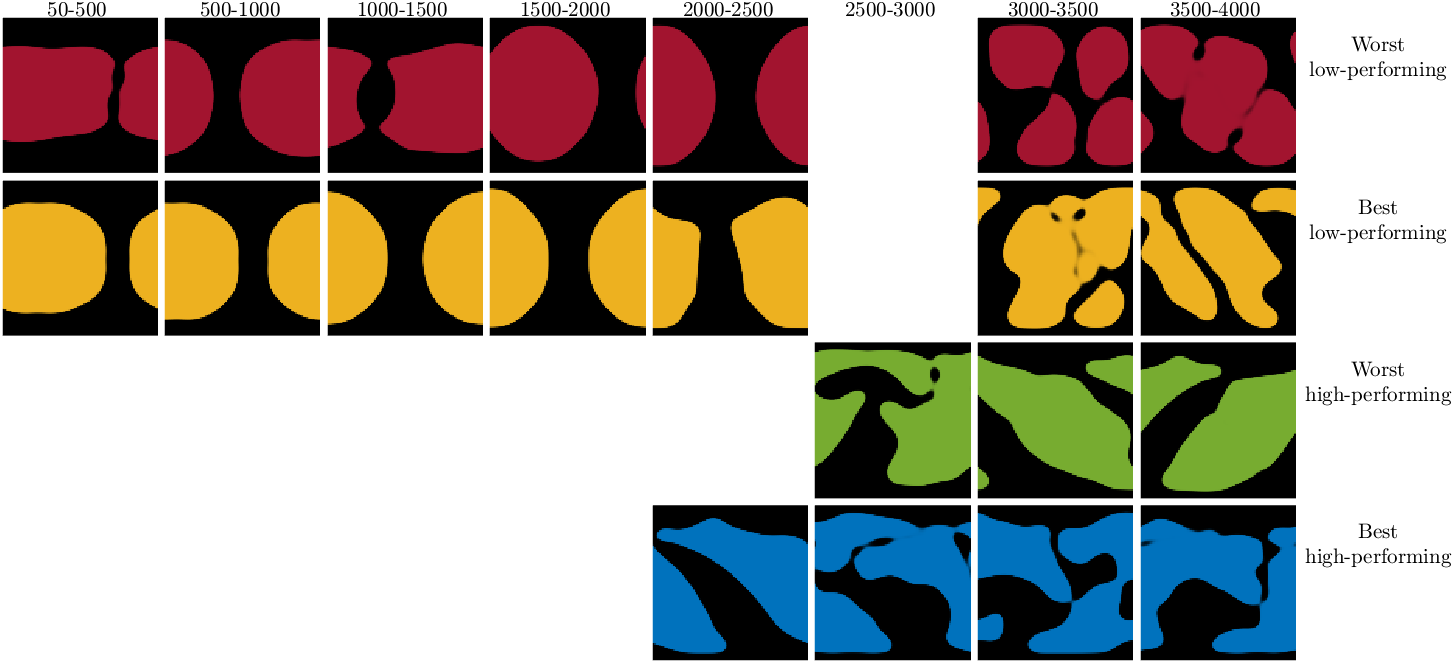}
\end{subfigure}\\
\begin{subfigure}[t]{\linewidth}
  \centering
  \captionsetup{width=\linewidth}
  \includegraphics[width=\linewidth]{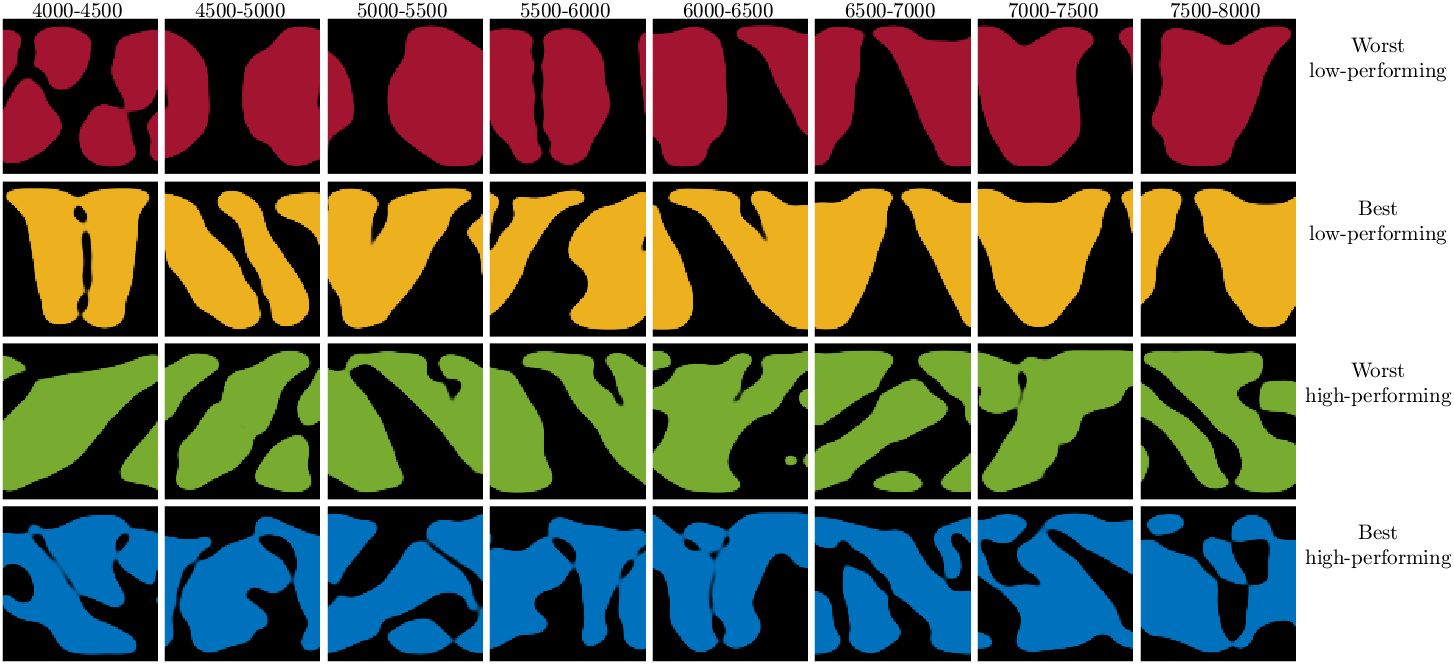}
\end{subfigure}%
\caption{Worst- and best-performing designs of the low- and high-performing optima for the sandwich panel case.}\label{fig:designs_baseline}
\end{figure}

\Cref{fig:designs_baseline} shows a selection of the optimal designs.
For every $\Delta \omega$, it plots two low-performing designs: the worst one in red and the best one in yellow.\ 
Also, for the high-performing optima, it plots the worst one in green and the best one in blue.\ 
A closer look at the low-performing designs below $2500\,\mathrm{Hz}$ reveals that it is filled with variants on the same design, namely a rigid vertical connection between top and bottom plate.\
Similarly, many designs for $\Delta\omega$ ranges between $2000\,\mathrm{Hz}$ and $5000\,\mathrm{Hz}$ consist of an askew connection between the top and bottom plate.\ 
For $\Delta \omega=[2000\,\mathrm{Hz},2500\,\mathrm{Hz}]$, this design is the (only) high-performing design.\
For $\Delta \omega$ ranges in $[3000\,\mathrm{Hz},4500\,\mathrm{Hz}]$, the design appears as the worst high-performing design (i.e., in green). 
Together with \Cref{fig:results_baseline}, this substantiates the existence of families of local optima whose performances both exceed and trail the mass law, depending on the targeted frequency range.\ 
The examples show that remaining trapped in low-performing families is detrimental.\ 
Hence, finding strategies to escape these regions is crucial.\

\section{Investigation of strategies for finding a high-performing optimum}\label{sec:strategy_investigation}

This section investigates three strategies (see \Cref{tab:strategy-overview}) for improving the chance of obtaining a high-performing local optimum.\ 
Of the two cases, the sandwich panel case is found to be the most difficult, exhibiting both a wide low-performing zone at low frequencies and many high-frequency transition regions.\
For this reason, it is chosen as the benchmark case against which the strategies are compared.\\

The strategies are judged against the following four criteria, in decreasing order of importance.

\begin{enumerate}
\item The \textbf{chance at a high-performing optimum}, approximated by $\tilde{P}_{\text{HP}}$.
\item \textbf{The average} $\mathbf{STL}$ \textbf{of the high-performing optima}, $\mathrm{STL}_{\text{HP}}=\sum_{i=1}^{N_{\text{HP}}} STL_i(\Delta \omega)/N_{\text{HP}}$.\ 
This captures the goal that high-performing optima should be as high-performing as possible.\
\item The \textbf{bias} towards a class of optima.\ 
This criterion penalizes strategies that always produces the same local optimum and is judged in a qualitatively way.
\item The \textbf{computation time}. Since the proposed strategies do not increase the time per iteration, this criterion is quantified by the average number of optimization iterations for the high-performing optima $N_{\text{iter}}^{\text{HP}}=\sum_{i=1}^{N_{\text{HP}}} N_{\text{iter}, i} / N_{\text{HP}}$.\ 
\end{enumerate}
Note that the second and fourth criteria are undefined for frequency ranges without high-performing optima ($N_{\text{HP}}=0$).\\

This section is built up as follows.\ 
\Cref{sec:gmin_comparison} studies the effect of exclusion methods that drive low-quality local optima towards infeasibility.\ 
For the application in this study, an upper bound on the stiffness is introduced to force the optimizer to produce compliant structures.\ 
\Cref{sec:fcont_comparison} investigates frequency shifts: by starting from a high targeted frequency range and decreasing it slowly, the greater density of high-performing optima at high frequencies is transplanted to low frequencies.\ 
\Cref{sec:robustness_comparison} examines a set of relaxation strategies that delay the required robustness against design changes.\ 
Finally, \Cref{subsec:comparison_strategies} provides a comparison of the studied strategies, highlighting advantages and disadvantages along all four considered criteria, and \Cref{subsec:generality_and_guidelines} provides a discussion on applicability to other problems.

\begin{table}[h!]
\caption{Overview of studied strategy variants.}
\label{tab:strategy-overview}
\centering
\begin{tabular}{m{0.25\linewidth}m{0.55\linewidth}m{0.2\linewidth}}
\toprule
\textbf{Strategy} & \textbf{Varied parameter(s)} & \textbf{Section} \\
\midrule
Exclusion & $J^{\min}_{\mathrm{conn}}$, a lower bound on the connectivity $J_{\mathrm{conn}}$ &  \Cref{sec:gmin_comparison} \\[2ex]
\midrule
Frequency shift &  Frequency range $\Delta \omega$ & \Cref{sec:fcont_comparison}  \\[2ex]
\midrule
Relaxation & Formulation (minmax or $1$-norm minimization), $\eta_e$, $\eta_d$ & \Cref{sec:robustness_comparison}  \\[2ex]
\bottomrule
\end{tabular}
\end{table}

\subsection{Lower bounds on the connectivity constraint} \label{sec:gmin_comparison}

This section explores a lower bound, $J^{\min}_{\mathrm{conn}}$, on the connectivity constraint to ensure (see \Cref{sec:constraints}) that $-1 < J^{\min}_{\mathrm{conn}} \leq J_{\mathrm{conn}} \leq 0$.\
Physically, this means that the optimal structure is \emph{required} to be less stiff against a self-weight load.\
The motivation for this strategy is the observation, valid for both the sandwich panel and cantilever case, that high-performing optima are less statically rigid than low-performing optima.\
By requiring less rigidity, low-performing optima become infeasible.\

Arguably, this violates the bias goal: by restricting the connectivity, stiff solutions are excluded and the optima are biased to be (more) compliant. If stiff, high-performing optima exist, this strategy can never find them.\\
Three variants of the proposed exclusion (E) strategy are considered.\
\begin{itemize}
    \item $\text{E}_1$: A constant, loose bound $J^{\min}_{\mathrm{conn}}=-0.5$.
    \item $\text{E}_2$: A constant, tight bound $J^{\min}_{\mathrm{conn}}=-0.05$.
    \item $\text{E}_3$: An adaptively increasing bound, starting at $J^{\min}_{\mathrm{conn}}=-1$ and going to $-0.05$, denoted as $J^{\min}_{\mathrm{conn}} \rightarrow -0.05$.
\end{itemize}
For the last variant, $\text{E}_3$, changes in $J^{\min}_{\mathrm{conn}}$ are triggered after $\beta_1 > 15$, allowing for spontaneous convergence to high-performing optima since these (see \Cref{sec:problem_illustration}) usually occur when $\beta_1 < 15$.\
Changes in $J^{\min}_{\mathrm{conn}}$ are applied by looking at the current value of $J_{\mathrm{conn}}$: $J^{\min}_{\mathrm{conn}}=\min_{n \in \mathbb{N}^+} -0.05n \quad \mathrm{s.t.}\quad J^{\min}_{\mathrm{conn}} > J_{\mathrm{conn}} $.\
If the optimizer finds an optimum with an active connectivity constraint ($J_{\mathrm{conn}} \approx 0$), $k=1$ and $J^{\min}_{\mathrm{conn}}$ is set immediately to $-0.05$.\\

\begin{figure}[h]
\centering
\begin{subfigure}[t]{0.5\textwidth}
  \centering
  \captionsetup{width=0.95\textwidth}
  \includegraphics[width=0.95\textwidth]{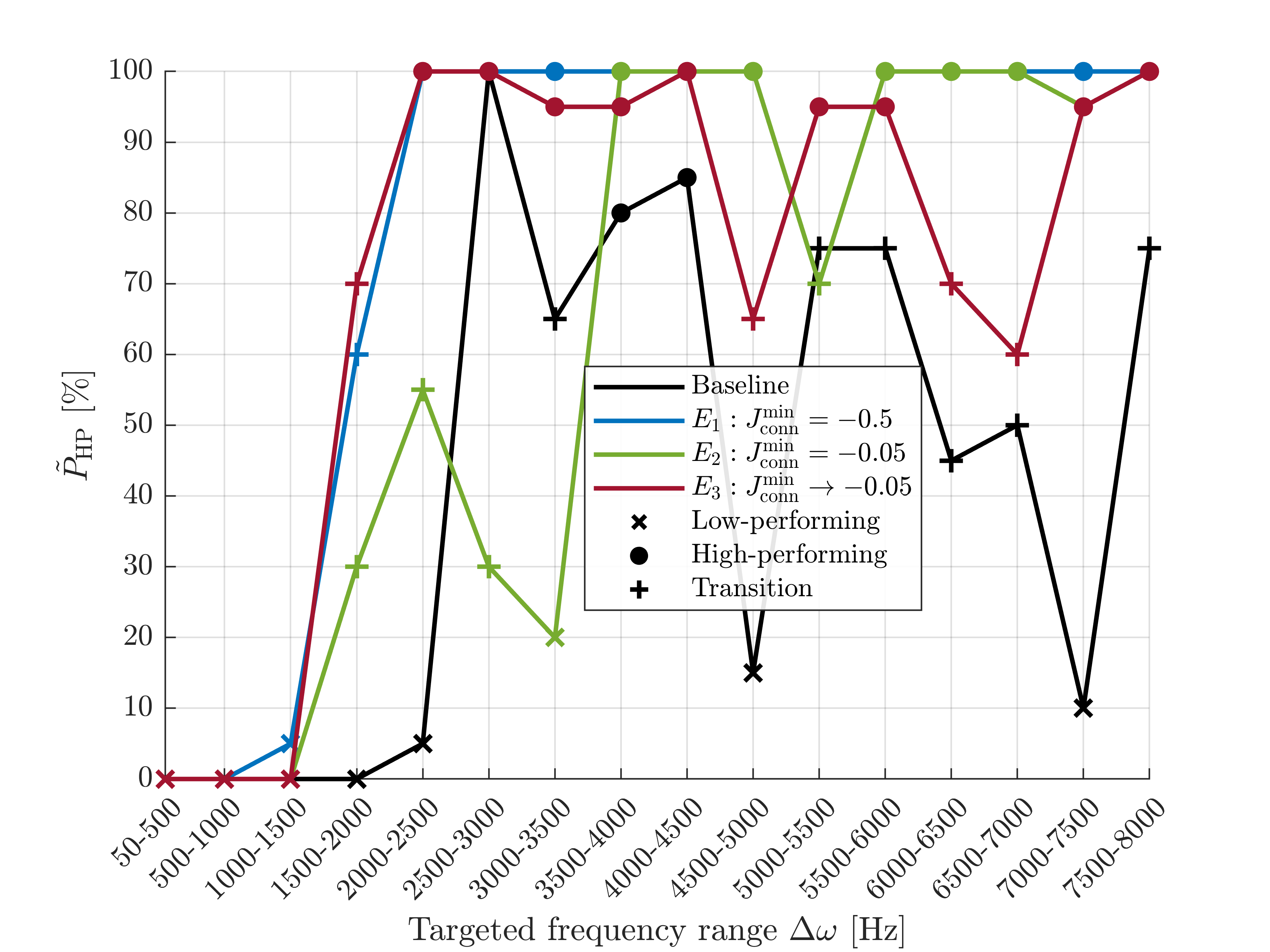}
  \caption{Chance of obtaining a high-performing optimum $\tilde{P}_{\text{HP}}$.}
  \label{fig:gmin_chance}
\end{subfigure}%
\begin{subfigure}[t]{0.5\textwidth}
  \centering
  \captionsetup{width=0.95\textwidth}
  \includegraphics[width=0.95\textwidth]{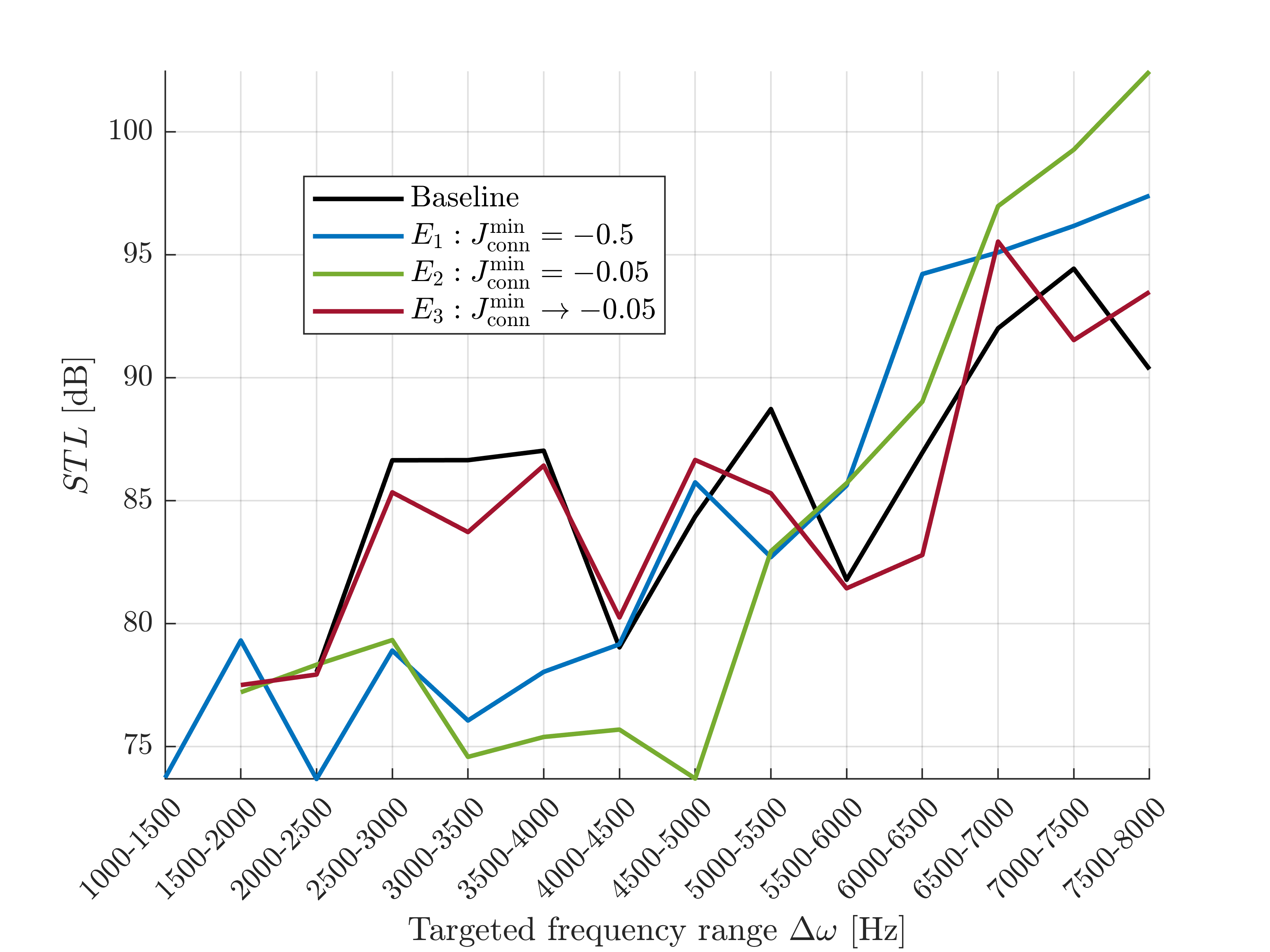}
  \caption{Average $\mathrm{STL}^*$ of the high-performing optima for all targeted frequency ranges.}
  \label{fig:gmin_avgSTLs}
\end{subfigure}%
\\
\begin{subfigure}[t]{0.6\textwidth}
  \centering
  \captionsetup{width=0.95\textwidth}
  \includegraphics[width=0.95\textwidth]{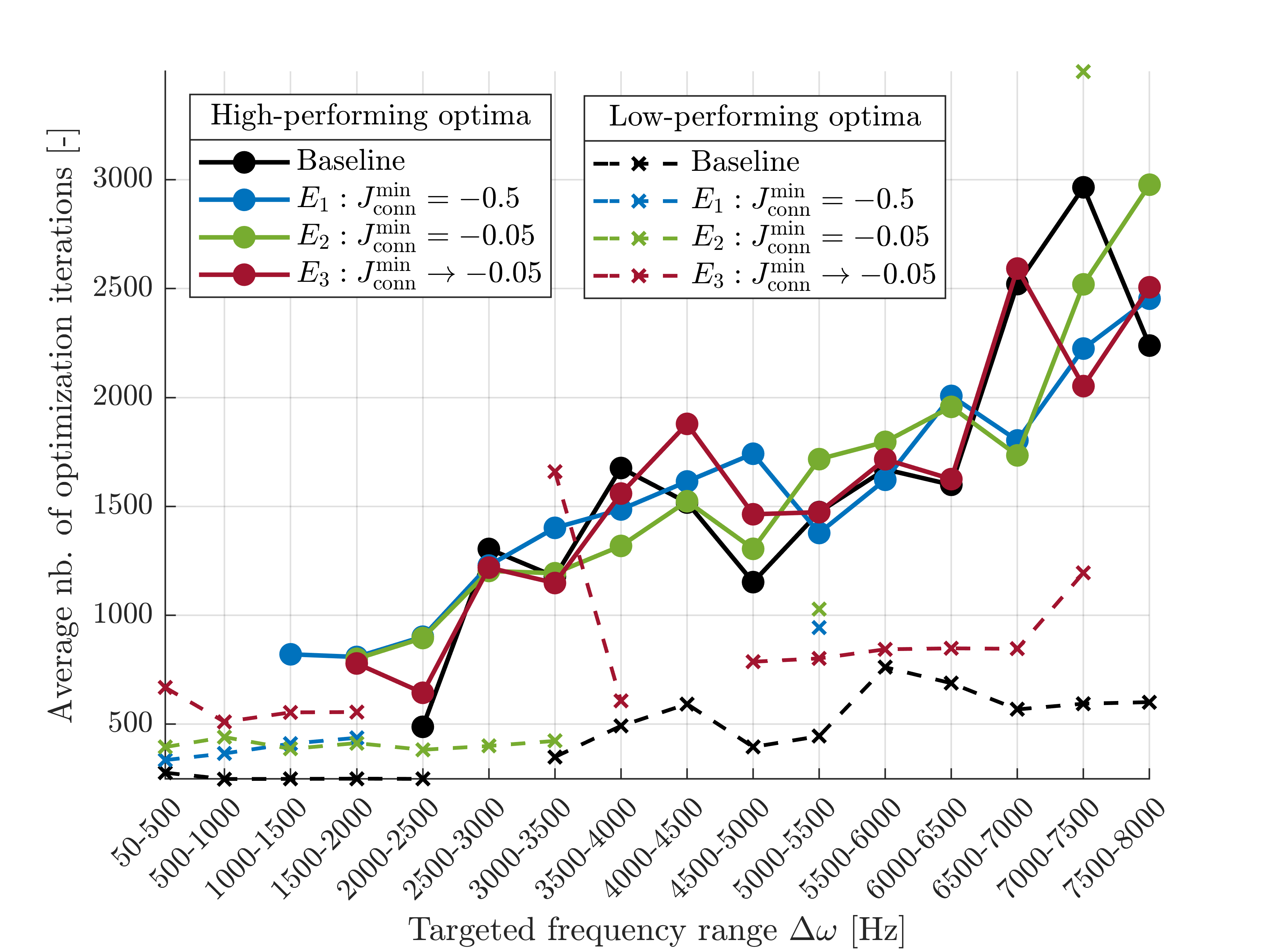}
  \caption{Average number of optimization iterations for the low- and high-performing optima.}
  \label{fig:gmin_niters}
\end{subfigure}%
\caption{Monte Carlo sampling results for the exclusion strategies, $E_1$ to $E_3$.}\label{fig:gmin_results}
\end{figure}

\Cref{fig:gmin_results} provides a comparison of the Monte Carlo sampling results for each variant.\
\Cref{fig:gmin_chance} visualizes $\tilde{P}_{\text{HP}}$ for the baseline and all three exclusion strategy variants.\ 
Variant $E_1$ both reduces the width of the low-performing low-performing zone by $1000\,\mathrm{Hz}$ and achieves $\tilde{P}_{\text{HP}}=100\%$ for all targeted frequency ranges above $2000\,\mathrm{Hz}$ except $[5000\,\mathrm{Hz}, 5500\,\mathrm{Hz}]$.\ 
Variant $E_2$ performs similar to $E_1$ above $3500\,\mathrm{Hz}$ but replaces the baseline's high-performing frequency range around $2500\,\mathrm{Hz}$ with a transition region.\ 
Variant $E_3$ does the opposite: it performs equally well as $E_1$ below $3500\,\mathrm{Hz}$ but has some transition regions above $3500\,\mathrm{Hz}$.\
Regarding the first criterion, exclusion variants can clearly outperform the baseline implementation, with variant $E_1$ outperforming $E_2$ and $E_3$.\\

The conclusion is more mixed when considering the second criterion, i.e., the $\mathrm{STL}^*$ of the high-performing optima.\ 
\Cref{fig:gmin_avgSTLs} shows that variant $E_3$ has comparable performance to the baseline: the average $\mathrm{STL}^*$ of the high-performing optima is not far from the baseline.\
In contrast, variants $E_1$ and $E_2$ differ significantly depending on the frequency range: around $3000\,\mathrm{Hz}$ the high-performing optima are $5-10\,\mathrm{dB}$ worse but at high frequencies they start outperforming the baseline by the same margins.\ 

Regarding the third criterion, i.e., bias, every exclusion variant a priori has a strict bias away from stiff designs.\ 
Between the variants, $E_1$ has a looser bound on the maximum stiffness than $E_2$ and $E_3$ and less bias.\\

Regarding the fourth and final criterion, computation time, theoretically none of the strategies incur a much higher computational cost per iteration.\ 
Furthermore, \Cref{fig:gmin_niters} shows that, except for some outliers due to a low number of samples, the number of optimization iterations is always comparable to the baseline. \ 

\subsection{Frequency range continuation} \label{sec:fcont_comparison}
This section investigates the frequency shift strategy from Olhoff and Du~\cite{olhoff2016_incremental_frequency} by changing the targeted frequency range $\Delta \omega$ during the optimization process.\ 
The motivation is the low-performing, low-frequency region seen for both cases in \Cref{sec:problem_quantification}.\ 
A natural idea is to start at a higher frequency range $\Delta\omega^*=[\omega_{-}+\omega^*, \omega_{+}+\omega^*]$, with $\omega^* > 0\,\mathrm{Hz}$ the excess frequency, and slowly remove the excess after a high-performing optimum is found.\
If the high-performing, high-frequency optimum is not ``lost'' during this process, the larger $P_{\text{HP}}$ at high $\Delta \omega$ is transplanted to the low-frequency, low-performing region.\\ 

\begin{figure}[h]
\centering
\begin{subfigure}[t]{0.5\textwidth}
  \centering
  \captionsetup{width=0.95\textwidth}
  \includegraphics[width=0.95\textwidth]{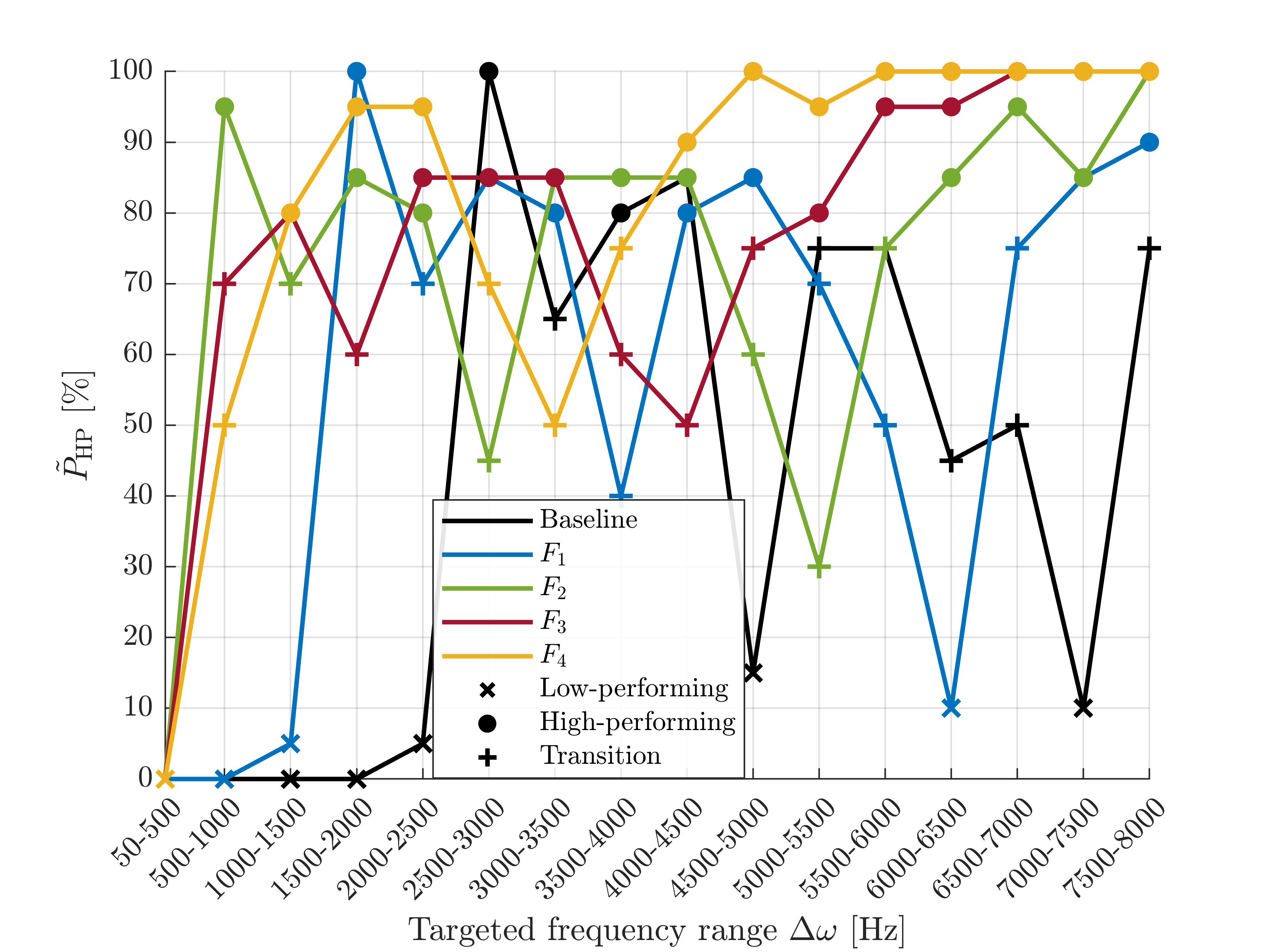}
  \caption{Chance of achieving a high-performing optimum $\tilde{P}_{\text{HP}}$.}
  \label{fig:fcont_chance}
\end{subfigure}%
\begin{subfigure}[t]{0.5\textwidth}
  \centering
  \captionsetup{width=0.95\textwidth}
  \includegraphics[width=0.95\textwidth]{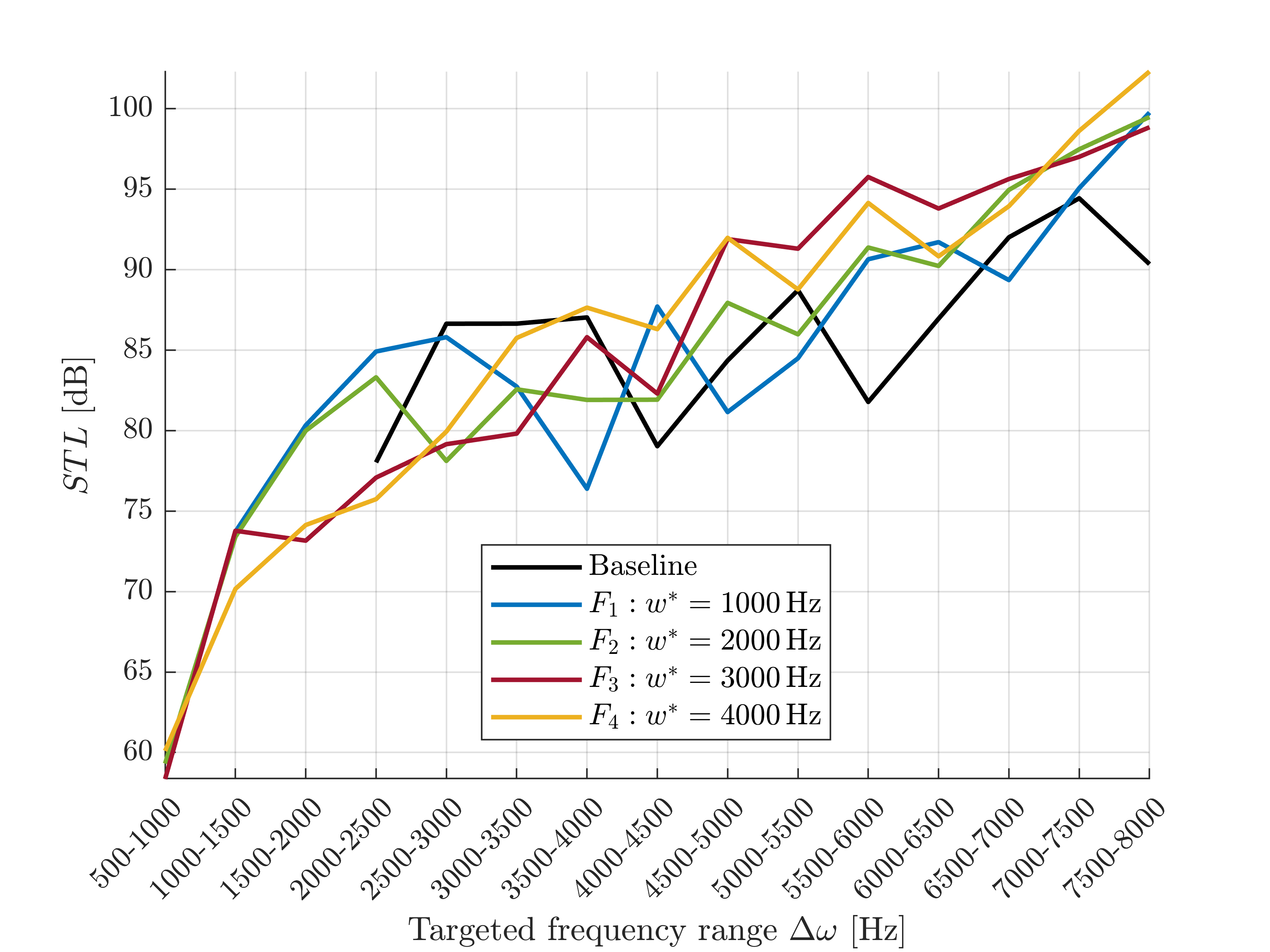}
  \caption{Average $\mathrm{STL}^*$ for high-performing optima.}
  \label{subfig:fcont_avgSTLs}
\end{subfigure}%
\\
\begin{subfigure}[t]{0.5\textwidth}
  \centering
  \captionsetup{width=0.95\textwidth}
  \includegraphics[width=0.95\textwidth]{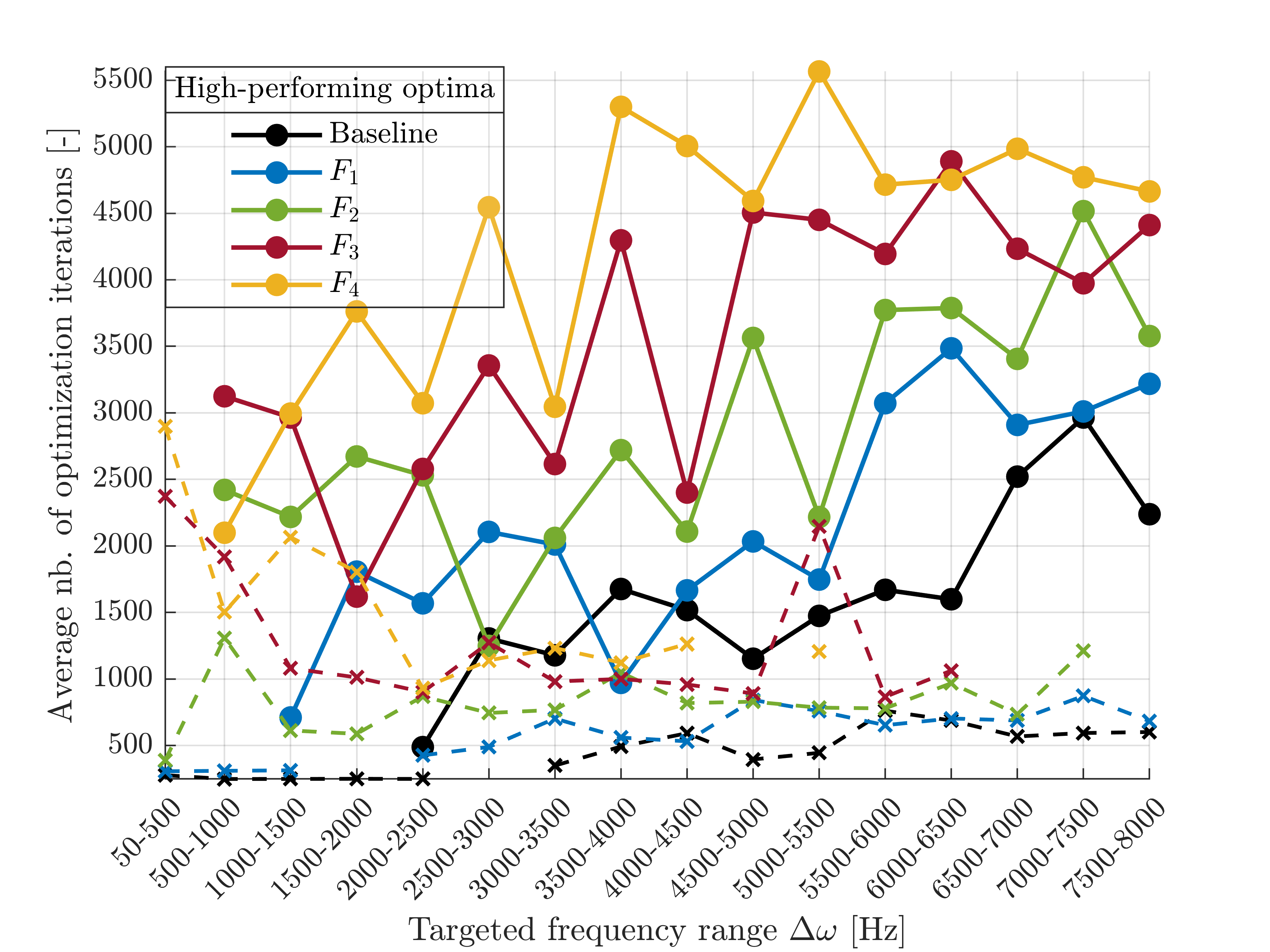}
  \caption{Average number of optimization iterations for the low- and high-performing optima. Low-performing optima curves are dashed.}
  \label{subfig:fcont_niters}
\end{subfigure}%
\caption{Monte Carlo sampling results for the frequency shift strategies, $\text{F}_1$-$\text{F}_4$.}\label{fig:fcont_results}
\end{figure}

Four variants of the frequency shift (F) strategy are tested: $\text{F}_n: \omega^*=1000n\,\mathrm{Hz}, n \in [1,\ldots, 4]$.\ 
The baseline routine starts with $\beta$ continuation until an intermediately converged, high-performing design is found: $\mathrm{STL}^*(\Delta\omega+\omega^*) \geq 1.15 \times STL_{ml}(\Delta\omega+\omega^*)$.\ 
Note that the criterion of $15\%$ above the mass law is a tunable hyperparameter that depends on the $STL$ separation between low- and high-performing classes of optima.\
This is further discussed in \Cref{subsec:generality_and_guidelines}.\
If this criterion is met, or if $\beta_1 > 50$, frequency continuation is initiated.\ 
The excess frequency $\omega^*$ is reduced by $100\,\mathrm{Hz}$ (or $1/5$th of the width of the considered frequency range) every time convergence is achieved.\ 
After $\omega^*=0$, $\beta$ continuation resumes.\
\Cref{fig:fcont_results} shows the results of the four frequency continuation strategies.

\Cref{fig:fcont_chance} allows to evaluate the first criterion, $\tilde{P}_{\text{HP}}$, for each frequency shift variant.\ 
The baseline's low-performing zone from $\Delta \omega = [50\,\mathrm{Hz},500\,\mathrm{Hz}]$ to $[2000\,\mathrm{Hz},2500\,\mathrm{Hz}]$ is decreased, depending on the size of the excess frequency $\omega^*$.\ 
That is, $F_1$ reduces it by its excess frequency of $1000\,\mathrm{Hz}$ and for variants $F_2$-$F_4$ only $\Delta \omega = [50\,\mathrm{Hz},500\,\mathrm{Hz}]$ remains low-performing.\ 
This reduction of the low-performing region is part of a larger shift of the baseline's $\tilde{P}_{\text{HP}}(\Delta \omega)$ curve.\ 
Namely, for frequency shift variant $F_i$, one sees that $\tilde{P}_{\text{HP}}^{F_i}(\Delta\omega) \approx \tilde{P}_{\text{HP}}^{\text{baseline}}(\Delta\omega+\omega^*)$.\ 
The exception is the lowest frequency range $\Delta \omega = [50\,\mathrm{Hz},500\,\mathrm{Hz}]$,  which implies a physical limit instead of a numerical one.\ 
Furthermore, if the excess frequency indeed shifts the baseline to the left, then \Cref{fig:fcont_chance} reveals the existence of a wide high-performing zone somewhere above $8000\,\mathrm{Hz}$, since $F_3$ and $F_4$ reach a plateau of $P_{\text{HP}}=100\%$ at high frequencies.\

Evaluation of the second criterion via \Cref{subfig:fcont_avgSTLs} shows that frequency shift methods generally outperform the baseline at higher frequencies (above $\pm 4000\,\mathrm{Hz}$).\ 
At low frequencies, variants $F_1$ and $F_2$ also see the leftward shift of the baseline's excellently performing designs ($\pm 85\,\mathrm{dB}$) between $2500$ and $4000\,\mathrm{Hz}$, causing them to outperform variants $F_3$ and $F_4$ below $3000\,\mathrm{Hz}$.\

Regarding computation time, measured via the number of optimization iterations shown in \Cref{subfig:fcont_niters}, a clear and marked increase with respect to the baseline is noticeable.\ 
This is due to the additional continuation steps: variant $F_n$ sees $10n$ frequency continuation steps on top of the usual $\beta$ continuation, each of which requires convergence.\ 
This explains why $F_4$ has the largest computation cost.\ 
The size of the increase is also noticeable: a doubling or even tripling of the computational cost is not unusual.\ 
Speedups could possibly be achieved through 1) more advanced $\omega^*$ reductions, e.g., based on the accompanying $\mathrm{STL}$ reduction, and 2) loosened convergence tolerances during frequency continuation steps.\
These are outside the scope of this work.\ 
Finally, for excessive $\Omega^*$, the high computational cost only leads to a physical limit: no high-performing optima were found for $\Delta \omega = [50\,\mathrm{Hz},500\,\mathrm{Hz}]$ but the number of iterations more than quadrupled for $F_3$ and $F_4$.\

\begin{figure}[h]
    \centering
    \includegraphics[width=0.95\textwidth]{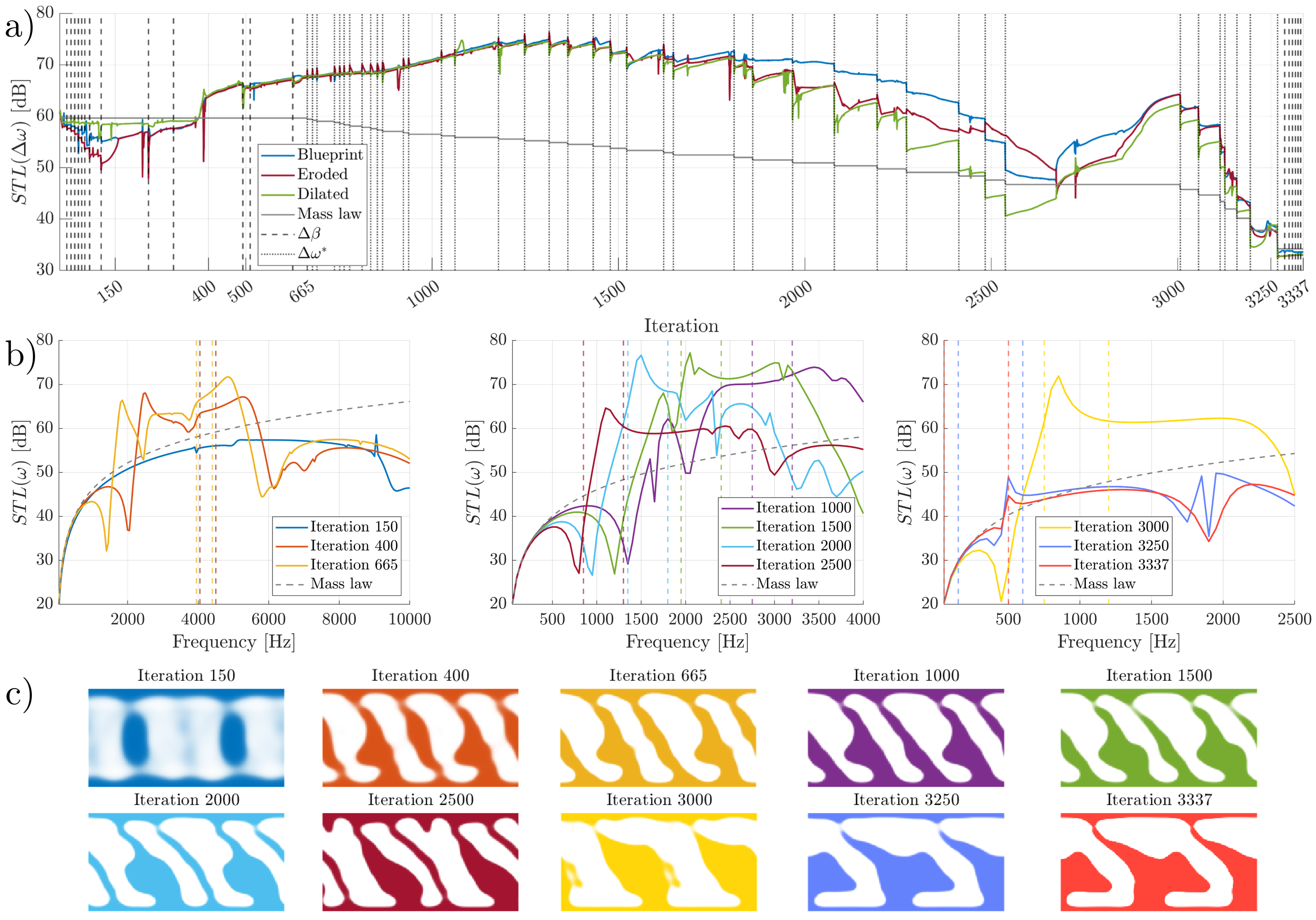}
    \caption{Optimization process overview of a run targeting $\Delta\omega=[50\,\mathrm{Hz},500\,\mathrm{Hz}]$ with frequency continuation ($\omega^*=4000\,\mathrm{Hz}$). a) Progress of $\mathrm{STL}_b(\Delta \omega)$, $\mathrm{STL}_e(\Delta \omega)$ and $\mathrm{STL}_d(\Delta \omega)$. b) $\mathrm{STL}_b(\omega)$ for selected blueprint designs. Vertical lines denote $\Delta \omega$ for the iteration with the same color. c) Selected blueprint unit cell layouts ($2\times$).}
    \label{fig:df4000_50Hz_run5}
\end{figure}

Regarding bias, frequency continuation strategies aim to draw high-frequency solutions towards lower frequencies and hence bias towards high-frequency optima.\ 
However, this is not guaranteed to be successful.\
\Cref{fig:df4000_50Hz_run5} shows an optimization run with variant $F_4$, targeting $\Delta\omega=[50\,\mathrm{Hz},500\,\mathrm{Hz}]$.\ 
Four stages can be discerned.\ 
In the first stage (iteration $1-665$), the usual optimization process is visible (cf. \Cref{fig:problem_illustration_dynamic}): stagnation for a few $\beta$ increases until a straight connection is formed (iteration $150$) and sharp improvement when this connection becomes slanted (iteration $400$).\ 
The $\mathrm{STL}$ improvement triggers the mass law criterion $\mathrm{STL}^* \geq 1.15 \times STL_{ml}$, which causes a frequency continuation step at iteration $665$.\ 
Here, the design contains a double askew connection with a decoupling frequency below $2000\,\mathrm{Hz}$.
The second stage (iteration $665-2500$) sees frequency continuation steps that decrease the frequency range from $ \Tilde 4000\,\mathrm{Hz}$ to around $1000\,\mathrm{Hz}$. The peak of the $\mathrm{STL}(\omega)$ curves move correspondingly to lie within this frequency band and the decoupling frequency is also pulled down to remain below the targeted frequency band.\ 
This keeps the $\mathrm{STL}(\Delta \omega)$ above the mass law until about iteration $2500$, where the decoupling frequency is only narrowly below $\Delta \omega$.\ 
The third stage (iteration $2500-3000$) combats this performance drop.\ 
The two slanted connections melt together into a connection with a tuned rotational inertia. This again pulls down the decoupling frequency and improves performance.\ 
However, stage four (iteration $3000-3337$) sees frequency continuations for which no corresponding decoupling frequency drop can be found.\ 
The performance again drops to the mass law and the routine terminates with several $\beta$ increases.\
Frequency shift strategies can thus bias towards suboptimal designs if a high-frequency solution cannot be adapted to a low-frequency target.\ 
More careful frequency continuation can combat the fallback to the mass law, but not the underlying bias towards solutions that exploit the decoupling frequency.\ 
In other words, the considered frequency shift variants are unlikely to produce high-performing designs with a decoupling frequency \emph{above} $\Delta \omega$.\\

Overall, frequency shift strategies can greatly reduce the low-performing zone at low frequencies.\ 
However, the right value for the excess frequency $\omega^*$ is not known a priori and a wrong choice can even decrease $P_{\text{HP}}$, such as when variant $F_2$ replaces the baseline's high-performing region at $\Delta \omega = [2500\,\mathrm{Hz}, 3000\,\mathrm{Hz}]$ with a transition region.\ 
This is doubly regrettable since frequency shifting incurs a significantly larger computational cost.\
In contrast, applying frequency continuation on the cantilever case would lead to much better results for two reasons.\
First, an obvious choice for $\omega^*$  exists, namely such that the initial frequency lies above the first natural frequency. 
Second, this choice guarantees one obtains a high-frequency high-performing solution with almost $100\%$ certainty (see \Cref{fig:BendsoeSigmund_baseline}).\
The only remaining challenge then is to remove the excess frequency without losing its high performance. 

\subsection{Design robustness and formulation} \label{sec:robustness_comparison}

\begin{figure}[h]
    \centering
    \includegraphics[width=0.95\linewidth]{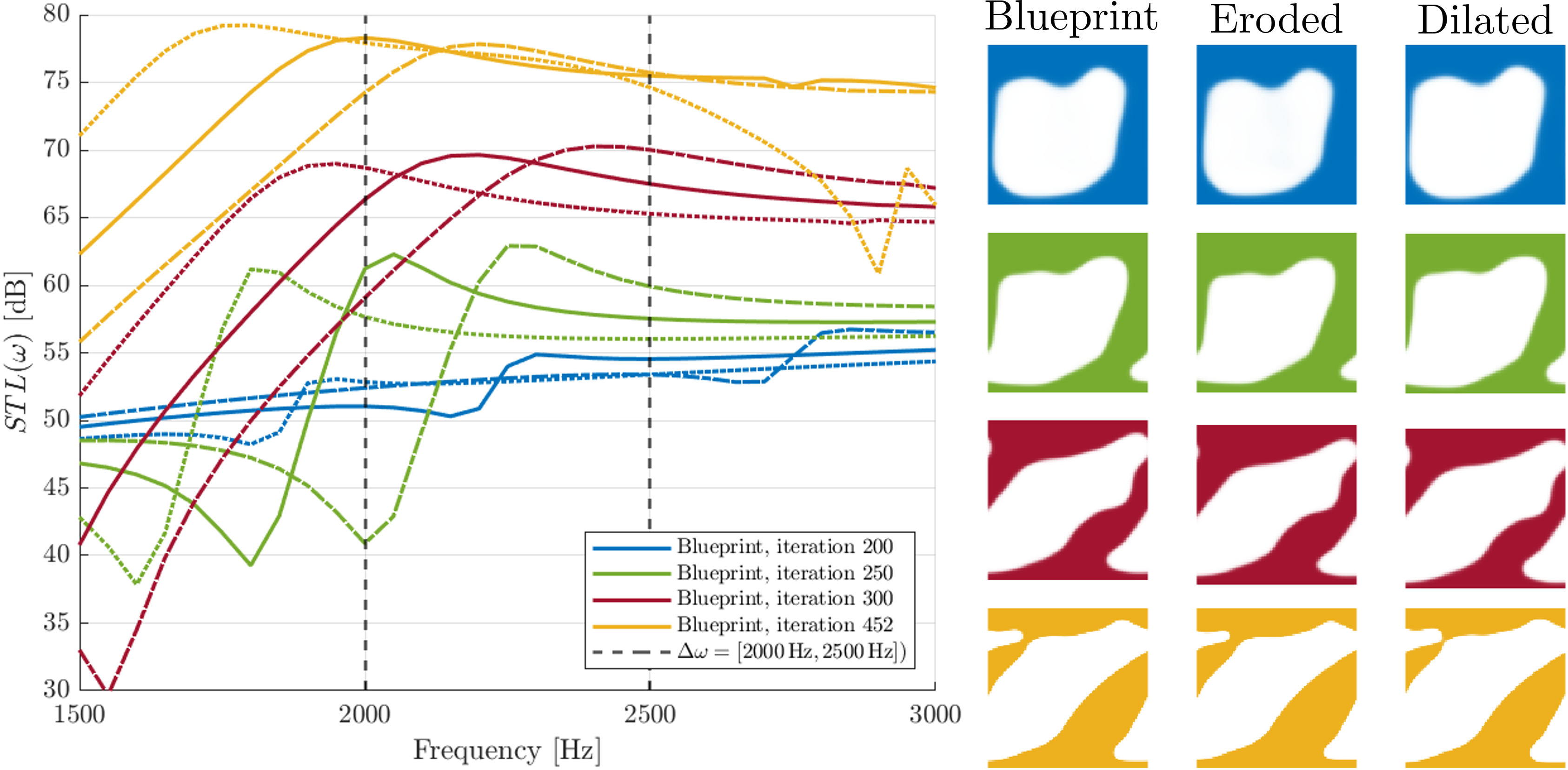}
    \caption{Eroded, blueprint and dilated designs (right) and corresponding $\mathrm{STL}(\omega)$ curves around $\Delta \omega$ for the run of \Cref{fig:problem_illustration_dynamic}. Eroded and dilated STL curves are dotted and dashdotted, respectively. }
    \label{fig:robustness_motivation}
\end{figure}

The third strategy tested in this work involves the relaxation of the design robustness.\ 
This is motivated by \Cref{fig:robustness_motivation}, which shows the $\mathrm{STL}(\omega)$ curves of the blueprint, eroded and dilated designs of some of the intermediate solutions of the optimization run of \Cref{fig:problem_illustration_dynamic}.\
All $\mathrm{STL}(\omega)$ curves show a trough followed by a peak.\ 
The eroded designs are shifted leftwards due to their lower mass and vice versa for the dilated designs.\ 
For iteration $200$, this causes the blueprint design to have a peak and trough in the middle of the targeted frequency range, whereas the eroded and dilated designs have their peaks and troughs left and right of that range, respectively.\ 
The design change from iteration $200$ to iteration $300$ sees an increase of $\mathrm{STL}(\Delta \omega)$ for all designs, but this is not self-evident because 1) the trough of the dilated design comes into the targeted frequency range while 2) (the tail end of) the peak of the eroded design moves out.\ 
Both effects cause a decrease of the performance of the worst design.\ 
However, the decoupling effect compensates for this and still causes an overall improvement.\ 

Based on the above, it is hypothesized that the minmax formulation's goal to improve the worst-performing design biases towards designs for which $\mathrm{STL}_b(\Delta\omega) \approx STL_e(\Delta\omega) \approx STL_d(\Delta\omega)$, which is largely true for mass-driven optima.\ 
Furthermore, a minmax formulation disallows large improvements of one design if accompanied by even a slight performance reduction of the worst design.\ 
Designs are prone to this phenomenon if the targeted frequency range lies between the STL decoupling peaks of the eroded and dilated designs, which might explain the baseline implementation's tendency to generate askew connections.\\

Two solutions are proposed.\ 
First, the minmax formulation can be replaced by an aggregated formulation, i.e., a formulation where the sum of the blueprint, eroded and dilated performance is optimized.\ 
\begin{linenomath}
\begin{equation}
\label{eq:aggregated_formulation}
\begin{aligned}
\min_{\bm{\xi}\in \mathbb{R}^{N_e}} \quad  J_b + J_e + J_d, \\
\end{aligned}
\end{equation}
\end{linenomath}
where the constraints are omitted for brevity.\ 
In this way, the aggregated formulation allows improvements in $\max STL_i(\Delta\omega), \, i \in [\mathrm{b, e, d}]$ as long as it compensates for the decrease in $\min STL_i(\Delta\omega), \, i \in [\mathrm{b, e, d}]$.\\ 

Second, the design robustness can be delayed.\ 
That is, only the blueprint design is considered at the start of the optimization and the eroded and dilated designs are introduced \emph{after} a high-performing optimum is found.\ 
This introduction of the eroded and dilated designs is done gradually by setting $\eta_{e}=\eta_b-\Delta \eta$, $\eta_{d}=\eta_b+\Delta \eta$ and slowly increasing $\Delta \eta$ from $0$ to its default value of $0.1$.\ 
The volume and connectivity constraint are still computed on the original dilated and eroded designs, respectively, to retain consistency with the baseline implementation.\ 

The two proposed solutions are now combined to form three different relaxation (R) strategies, $\text{R}_1$, $\text{R}_2$ and $\text{R}_3$.\

Strategy $\text{R}_1$ starts with the aggregated formulation, maximizing the sum of the performances of the blueprint, eroded and dilated designs.\ 
To ensure consistency with the baseline implementation, the minmax formulation is reintroduced after the last $\beta$ continuation step.\ 
Strategy $\text{R}_1$ is denoted by ``$\sum \rightarrow$ minmax''.\

Strategy $\text{R}_2$ starts by considering only the blueprint STL.\ 
Then, design robustness (with low $\Delta \eta=0.02$) with the aggregated formulation is introduced after the criterion of \Cref{sec:fcont_comparison} is met: a converged design with $\mathrm{STL}^* \geq 1.15 \times STL_{ml}$ or $\beta>50$.\ 
Next, the design robustness increases gradually: $\Delta \eta$ goes from $0.02$ to $0.1$ in four steps of $0.02$, each step requiring convergence (see \ref{app:continuation_convergence}).\ 
Finally, just like strategy $A$, the formulation switches from aggregated to minmax after the last $\beta$ increase.\
Strategy $\text{R}_2$ is denoted by ``$\sum \rightarrow \Delta \eta \uparrow \rightarrow $ minmax''.\

Strategy $\text{R}_3$ employs the same $\Delta \eta$ continuation strategy as strategy $\text{R}_2$ but uses the minmax approach for the whole optimization process.\ 
All three strategies end with the same design robustness and formulation as the baseline approach, allowing for a fair comparison.\ 
Furthermore, strategy $\text{R}_1$ differs from strategy $\text{R}_2$ only in the $\Delta \eta$ continuation and strategy $\text{R}_2$ differs from strategy $\text{R}_3$ only in the formulation.\ 
This allows for an unbiased conclusion as to the effect of 1) the design robustness and 2) the formulation on the chance at a high-performing optimum.\
Strategy $\text{R}_3$ is denoted by ``minmax $\rightarrow \Delta \eta \uparrow$''.

\begin{figure}[h]
\centering
\begin{subfigure}[t]{0.5\textwidth}
  \centering
  \captionsetup{width=0.95\textwidth}
  \includegraphics[width=0.95\textwidth]{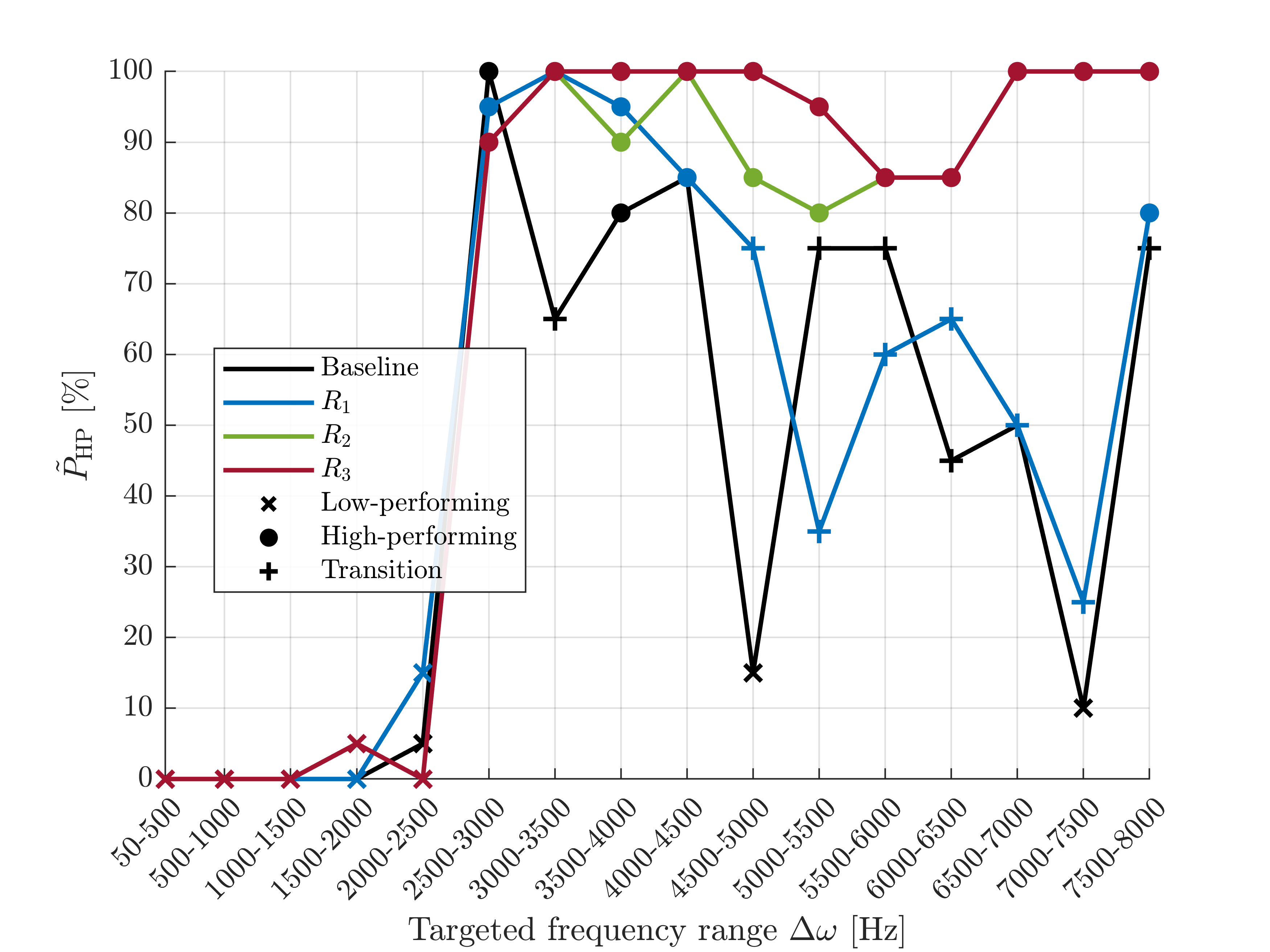}
  \caption{Chance of achieving a high-performing optimum $\tilde{P}_{\text{HP}}$ ($y$-axis) as a function of the optimized frequency range $\Delta\omega$ ($x$-axis).}
  \label{subfig:dynamic_chance_robustness}
\end{subfigure}%
\begin{subfigure}[t]{0.5\textwidth}
  \centering
  \captionsetup{width=0.95\textwidth}
  \includegraphics[width=0.95\textwidth]{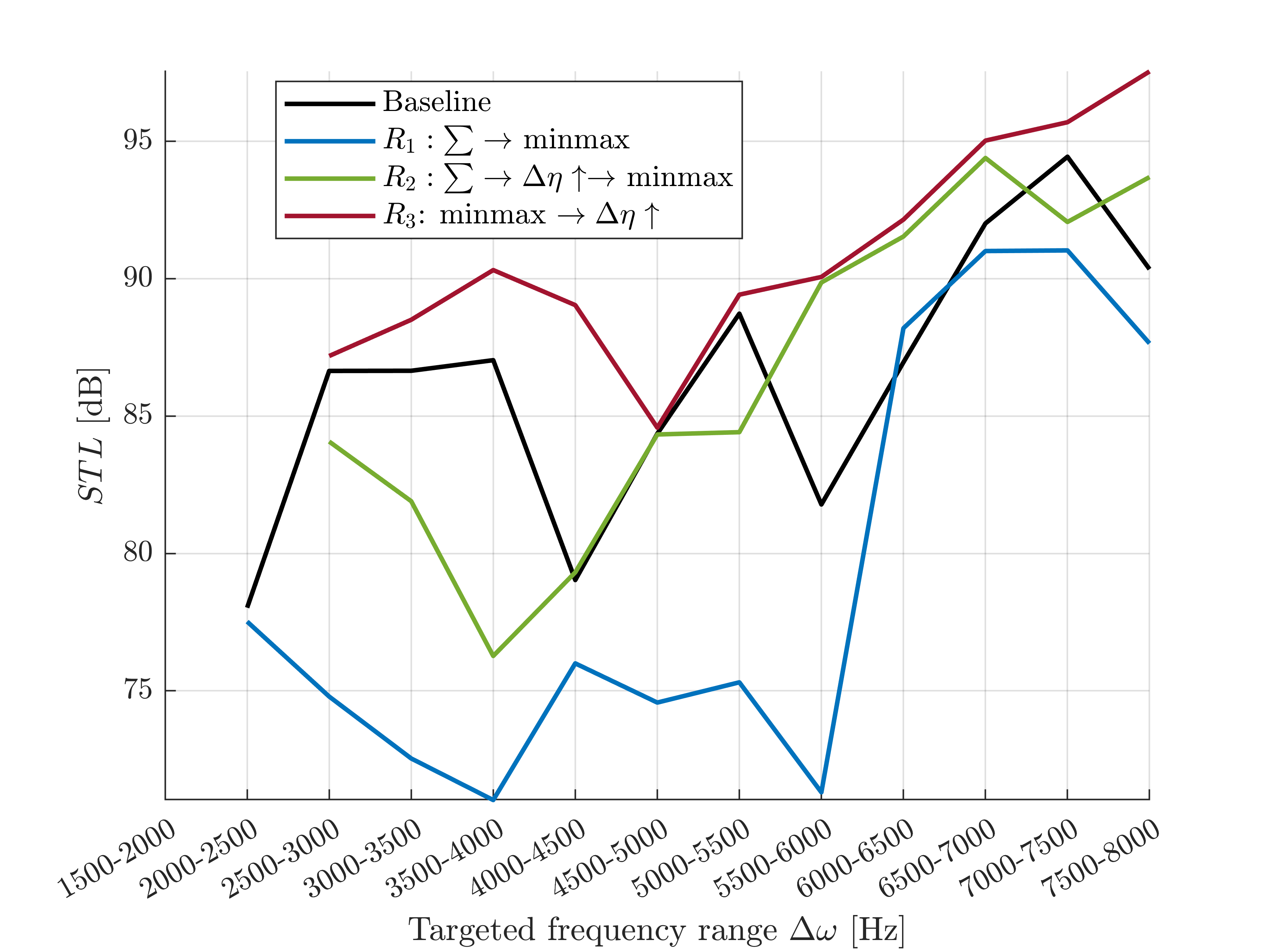}
  \caption{Average $\mathrm{STL}^*$ for low- and high-performing optima for all targeted frequency ranges.}
  \label{subfig:robustness_avg_STL}
\end{subfigure}\\
\begin{subfigure}[t]{0.5\textwidth}
  \centering
  \captionsetup{width=0.95\textwidth}
  \includegraphics[width=0.95\textwidth]{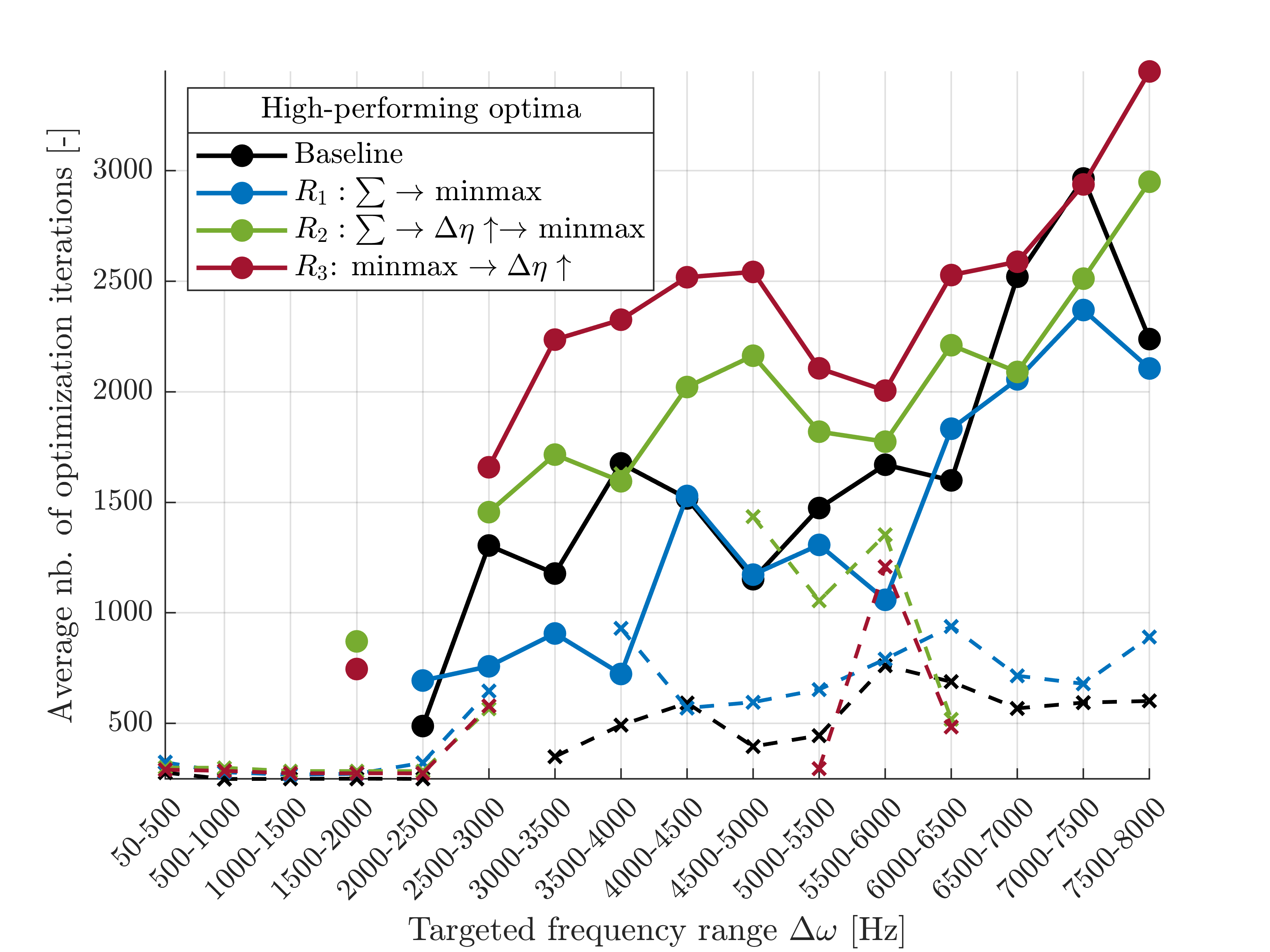}
  \caption{Average number of optimization iterations for the low- and high-performing optima.}
  \label{subfig:robustness_niters}
\end{subfigure}%
\caption{Monte Carlo sampling results for the relaxation strategy variants, $\text{R}_1$, $\text{R}_2$ and $\text{R}_3$.}\label{fig:results_robustness}
\end{figure}

\Cref{fig:results_robustness} shows the results of the three relaxation strategy variants.\ 
\Cref{subfig:dynamic_chance_robustness} shows that although none of the three variants are able to decrease the width of the low-frequency, low-performing zone, both $R_2$ and $R_3$ yield an unbroken high-performing region above $3000\,\mathrm{Hz}$.\ 
On the other hand, $R_1$ shows no $\tilde{P}_{\text{HP}}$ improvement.\ 

\Cref{subfig:robustness_avg_STL} shows that the performance of the high-performing optima is best for variant $R_3$, which always outperforms the baseline, albeit sometimes narrowly.\ 
Variant $R_2$ has, on average, comparable performance to the baseline: slightly worse at low frequencies and slightly better at higher frequencies.\ 
Finally, variant $R_1$ is often at least $5\,\mathrm{dB}$ worse than the baseline.\\ 

\Cref{subfig:robustness_niters} shows the cost of the relaxation strategy variants.\ 
Whereas $R_1$ generally requires less iterations than the baseline, $R_2$ requires slightly more and $R_3$ sometimes even sees more than a doubling of the number of iterations.\ 
This is caused by the larger number of continuation steps: every increase of $\eta$ requires convergence before another step is taken.\ 
It is worth noting that, compared to the baseline, the computational cost per iteration is lower for strategies that delay design robustness ($R_2$ and $R_3$) since only the blueprint design must be evaluated.\
However, due to parallelization, all evaluations occur simultaneously and hence no distinct reduction of the computation time per iteration is noticeable.\

Regarding bias, a study of the optimal designs reveals a wide array of design types and no clear bias.\ 
Possibly, the initial absence of design robustness of variants $R_2$ and $R_3$ allow them to better navigate towards high-performing designs.\ 

Overall, it is clear that relaxation strategies can strongly improve both the chance at a high-performing optimum and their performance.\
Of the three variants, $R_3$ comes out on top.\
This variant sets $\Delta \eta=0$ until a high-performing optimum is found.\
The associated computational cost increases due to the additional continuation steps, but only by at most a factor of two.\

\subsection{Comparison} \label{subsec:comparison_strategies}

The three strategy classes (exclusion, frequency shift and relaxation) are now compared for each of the four evaluation criteria. A summary is provided in \Cref{tab:strategy-comparison}.\\

For $P_{\text{HP}}$, only frequency shift techniques can significantly shrink the low-performing zone at low frequencies, until reaching the physical limit,  namely, the decoupling frequency.\ 
Exclusion strategies offer slight improvement, while relaxation methods have no effect.\ 
In contrast, for the transition regions and isolated low-performing regions, the trend reverses: a poor choices of $\omega^*$ degrades $P_{\text{HP}}$, whereas exclusion and relaxation strategies turn these regions into high-performing ones.\
\\

For $\mathrm{STL}_{\text{HP}}$, exclusion and frequency shift strategies have mixed effects: improvement at high frequencies but harming low ones. Relaxation strategy variant $R_3$ consistently improves $\mathrm{STL}$ across the entire frequency range.\\ 

Regarding bias, exclusion strategies favor compliant designs since stiff designs are a priori excluded. Frequency shift strategies bias toward high-frequency solutions tuned to the target frequency, while relaxation strategies only show a desirable bias away from the low-performing designs.\\

Finally, regarding computation cost, exclusion and relaxation strategies add little overhead, but frequency shifting can double or triple the number of iterations.\
Careful convergence and continuation management is needed to mitigate this.

\begin{table}[h!]
\caption{Strategy comparison across the chosen evaluation criteria.}
\label{tab:strategy-comparison}
\centering
\begin{tabular}{m{0.1\textwidth}m{0.3\textwidth}m{0.2\textwidth}m{0.12\textwidth}m{0.14\textwidth}}
\toprule
\textbf{Strategy} & \multicolumn{4}{c}{\textbf{Criteria}} \\ \\
\cmidrule(lr){2-5}
& $P_{\text{HP}}$ & $\mathrm{STL}_{\text{HP}}$ & Bias & Computational cost \\
\midrule
\\
Exclusion & \RaggedRight{Turns transition regions into high-performing regions, reduces low-performing, low-$\Delta \omega$ zone.} & \RaggedRight{Improvement at higher $\Delta \omega$, decrease at lower $\Delta \omega$.} & \RaggedRight{Towards compliant designs.} & \RaggedRight{No difference.}\\[4ex]
\midrule
Frequency shift & \RaggedRight{Reduces low-performing, low-$\Delta \omega$ zone but can replace high-performing regions with transition regions for poor $\omega^*$ choice.} & \RaggedRight{Improvement at higher $\Delta \omega$, decrease at lower $\Delta \omega$.} & \RaggedRight{Towards optima of higher $\Delta \omega$.} & \RaggedRight{Often doubles or triples cost, depending on setup.} \\[4ex]
\midrule
Relaxation & \RaggedRight{Turns transition regions into high-performing regions. No influence on low-performing, low-$\Delta \omega$ zone.} & \RaggedRight{Can improve performance across all frequencies (variant $R_3$).} & \RaggedRight{/} & \RaggedRight{Modest increase (usually $< 2\times$).} \\[4ex]
\bottomrule
\end{tabular}
\end{table}

\subsection{Applicability and guidelines} \label{subsec:generality_and_guidelines}

Due to the nonconvex nature of topology optimization problems, optima are highly dependent on the underlying physics, the optimization problem, the initial guess, associated hyperparameters, ...\
Nonetheless, the proposed strategies could be implemented and adapted for contexts other than the one discussed in this study.\ 
To this end, the following paragraphs provide a discussion and guidelines on applicability to other contexts.\

\subsubsection{Influence of hyperparameters}
In topology optimization, physical problem definition parameters such as element and unit cell size, material properties and  frequency range, ... heavily influence the available optima.\ 
For example, high-performing sandwich panels rely on decoupling and resonance effects, which are frequency dependent and hence depend in turn on size, mass and stiffness.\
To minimize the dependence of the frequency, the design domain size and the material properties are fixed and the selected frequency ranges are carefully chosen to ensure ranges both below and above the initial structure's first natural frequency.\ 
Furthermore, the minimum length scale is small enough in comparison to the unit cell length to allow for sufficient design freedom and the minimal stiffness (governed by $\mu_{sw}$) was chosen low enough to ensure that there was room to shift natural frequencies downwards.\ 
Based on the above, changes to the material, domain size or target frequency will likely lead to optima that rely on the same physical effects (resonance and decoupling) but occur at different frequencies, without fundamentally altering the conclusions of the proposed strategies.\

Four hyperparameters, related to the continuation and convergence criteria, require further motivation.\ 
First, tests prior to this study generally show that the convergence criterion is relatively unimportant: continuation should occur whenever the design seems converged, but $P_{\text{HP}}$ is not sensitive to the exact convergence hyperparameters.\ 
Second, for the sandwich panel case, the MMA move limit strongly influences the speed of convergence and numerical stability: too high and the routine becomes oscillatory, too low and convergence is slow.\ 
The adaptive MMA move limit scheme discussed in \ref{app:TO} is used to automatically tune this during the optimization.\
Third, some strategies rely on the criterion  $\mathrm{STL}^* \geq \alpha \times STL_{ml}$ to check if a design is high-performing.\ 
As long as low- and high-performing $\mathrm{STL}$ values are sharply separated (as in \Cref{fig:results_baseline}), the conclusions of this study are not sensitive on the exact value of $\alpha$.\
Fourth, frequency shift and relaxation methods rely on a ``step size'' parameter, which determines speed with which solutions are transplanted from, e.g., high to low frequency.\
To minimize the number of high-performing ``lost along the way'', small step sizes are used: the incremental frequency shift is only $1/5$th of the width of the frequency range and the design robustness is toggled from $\eta=0$ to $\eta=0.1$ in five steps of $\Delta \eta=0.02$.\
These step sizes could be further improved, e.g., by an adaptive scheme that tracks the chance in objectives and constraints, but suffice for the scope of this study.\

Overall, the convergence and continuation parameters in this work either do not have a large influence or are adapted automatically (MMA move limit) or are chosen conservatively (e.g. frequency step size).\

\subsubsection{Applicability to other problems}

For the cantilever and sandwich panel case, the appearance of eigenfrequencies endow the objective function space (dynamic compliance and STL, respectively) with a combination of flat regions with poor local optima and sparse, sharp valleys with high-performing optima.\  
Although the Monte Carlo study was only applied to those two cases, it is likely that this occurs in different dynamic contexts as well.\
Examples include:
\begin{itemize}
\item topology-optimized resonator design, such as for metamaterial resonators~\cite{lirias4249142} or MEMS~\cite{GIANNINI2022104352}, since resonators are sharp local optima;
\item work on bandgap maximization from, e.g. Wu et al.~\cite{Wu2025BandGap} and Cool et al.~\cite{COOL2025117744}, as these are prone to local minima and also use the design robustness scheme from Wang et al.~\cite{wang2011projection};
\item STL maximization in other contexts, e.g., for composite sandwich structures on the part-scale~\cite{HU2025118109}, due their similar optimization problem structure;
\item  works on frequency response matching, such as the eight objective functions considered by Zhao et al.~\cite{ZHAO2019700}, since these often lead to structures with internal resonators.
\end{itemize}

For each of the above examples, it could be of interest to apply the proposed strategies.\ 
However, adaptation to other contexts requires careful selection of hyperparameters and changes to the formulation.\
To facilitate adaptation, a set of lessons learned is offered, which can serve as guidelines and inspiration for researchers faced with similar challenges:\

\begin{itemize}
\item Before attempting to increase $P_{\text{HP}}$, it is crucial to confirm that high-performing optima exist.\
This can done through physical insight, e.g., by simplifying the problem, and via a limited parameter selection study to investigate, e.g.,  the loading frequency (low vs high) and the design freedom (filter radius, element size, ...).
\item Frequency shift methods can work very well to transplant high-frequency solutions to low frequency ranges.\
A limited Monte Carlo study (i.e. $N=3$ samples for frequencies below, on and above the first resonance) can (i) show if this is the case in another context and (ii) provide a good choice for the excess frequency $\omega^*$. \
\item If the class of poor local optima is well-understood, an exclusion strategy is viable.\
For example, if poor optima are stiff, as in this study, an upper stiffness bound could be used.\
Other contexts require tailored constraints.\
However, note that exclusion strategies are shown to work best for transition regions, i.e., to increase $P_{\text{HP}}$ in regions where $P_{\text{HP}}$ is already nonzero.
\item Relaxation strategies are useful for strongly constrained optimization problems.\
Examples include minmax formulations, like in this study, and stress-constrained problems due to the numerous local constraints.\ 
\item If an analytical expression for poor performance (like the mass law in this study) is not available, one could alternatively use an average of the performance in the first few iterations.
\item Continuation step sizes should initially be chosen small and increased gradually until high-performing optima are lost.\ 
If computation cost is critical, an adaptive step size can be considered.
\end{itemize}

Overall, the above guidelines and lessons learned can serve as a useful starting point for adaptation to a wide range of dynamic optimization contexts.

\section{Conclusions} \label{sec:conclusions}

This work divides methods that aim to improve $P_{\text{HP}}$, the probability of finding high-performing optima in dynamic topology optimization, into four categories: (i) global optimization, (ii) exclusion, (iii) relaxation and (iv) frequency shift methods.\
The work considers two case studies.\
The first is the dynamic compliance minimization of the reinforced cantilever of Bendsøe and Sigmund, for which high-performing optima are only found at sufficiently high driving frequencies.\
The second case considers the $\mathrm{STL}$ maximization of a sandwich panel.\
A Monte Carlo sampling methodology identifies two distinct subproblems: (i) a zone at low frequencies completely absent of high-performing optima ($P_{\text{HP}}=0\%$), similar to what's known for the cantilever case, and (ii) a high-frequency zone where $P_{\text{HP}}$ is often below $50\%$.\
This is caused by classes of numerically attractive yet low-performing optima. \\

To improve upon the Monte Carlo results, which is a global optimization strategy, the study presents strategies to maximize $P_{\text{HP}}$: a novel exclusion and a novel relaxation variant, as well as several instances of the frequency shifting technique from Olhoff and Du~\cite{olhoff2016_incremental_frequency}.\ 
For each strategy, various instances are developed, tested and compared for the sandwich panel case.\

The first strategy, exclusion, adds constraints to render the class of low-performing optima infeasible.\ 
Since low-performing optima maximize the sandwich panel stiffness, an upper bound is placed on stiffness.\ 
The second strategy, frequency shifting, optimizes first at higher frequencies and then gradually reduces the excess frequency, thereby transplanting high-performing optima into regions that lack such designs.\
The third strategy, relaxation, first solves a ``simpler'', less constrained formulation and uses the result as a starting point for gradually more constrained problem variations.\ 
Since $\mathrm{STL}$ maximization is strongly limited by the need for robustness against design changes, this requirement is initially relaxed.\\

On top of (i) their effect on $P_{\text{HP}}$, strategies are also compared by their ability to (ii) raise the $\mathrm{STL}$ of the high-performing optima, (iii) avoid bias and (iv) maintain reasonable computation cost. The main conclusion is that frequency shift techniques are particularly effective for increasing $P_{\text{HP}}$ in the low-performing zone at low frequencies, but incur up to a tripling of the computational cost and introduce a bias towards high-frequency designs. Conversely, exclusion and relaxation techniques can increase $P_{\text{HP}}$ to $100\%$ in the mid- to high-frequency zone for only a limited increase in computation time.\\

Naturally, the proposed strategies and their hyperparameters are tailored to the studied cases, which are linear and 2D.\
Extensions to 3D or non-linear models could require adaptations or invalidate the effect of the continuation methods.\
Nevertheless, the core ideas (exclusion, relaxation and frequency shifting) are applicable to many other settings and the continuation strategies are simple enough to be readily adapted.\ 
To this end, the provided guidelines offer researchers a practical basis to identify suitable methods and adapt them to their application.\ 
Ultimately, the proposed strategies could enhance the topology optimization community's ability to obtain high-performing designs without excessive amounts of reruns.\

\section*{Statements and Declarations}
\subsection*{Funding}
Resources and services were provided by the VSC (Flemish Supercomputer Center), funded by the Research Foundation - Flanders (FWO) and the Flemish Government.\
The research of V.\ Cool (fellowship no.\ 1213925N) is funded by a grant from the Research Foundation - Flanders (FWO).\
Internal Funds KU Leuven are gratefully acknowledged for their support.\ 
The Research Foundation - Flanders (FWO) is gratefully acknowledged for its support through research grant no. G0F9922N.

\subsection*{Conflict of Interest}
On behalf of all authors, the corresponding author states that there is no conflict of interest.

\subsection*{Author Contributions}
All authors contributed to the study conceptualization and design. Data collection, analysis and writing was performed by the corresponding author. All authors commented on previous versions of the manuscript. All authors read and approved the final manuscript.

\subsection*{Ethics approval and Consent to participate}

Not applicable.

\subsection*{Data Availability and Replication of Results}

The details of the proposed strategies are described to facilitate the replication of the results.\
 Data will be made available on request.

\appendix
\section{Topology optimization framework}
\label{app:TO}
This appendix provides further numerical details regarding the optimization framework for the two cases solved in this work.\

\subsection*{Case specification}

The considered cantilever has $L=1\,\mathrm{m}$ and is meshed with $200 \times 200$ equisized, bilinear finite elements.\ 
The outer $b=10$ elements are passive elements kept solid during optimization and the applied force has size $F=100\,\mathrm{N}$ and is spread out over a surface of $10$ elements.\ 
 The solid material has Young's modulus $E=5\,\mathrm{GPa}$, Poisson's ratio $\nu=0.3$ and density $\rho=1000\,\mathrm{kg/m^3}$.\
The void material has $E=5\,\mathrm{kPa}$ and $\rho=1.e-9\,\mathrm{kg/m^3}$ to prevent spurious low-density modes~\cite{Pedersen2000}.\
The structural damping coefficient is $\alpha=0.1$.\
The framework uses RAMP interpolation with $q=3$ for the Young's modulus and linear interpolation for the density.\\

The sandwich panel unit cell has dimensions $L_x \times L_y=50\,\mathrm{mm} \times 50\,\mathrm{mm}$, with a plate thickness of $5$~mm.\
A rectangular $N_x \times N_y$ Finite Element grid is employed for a total of $N_e$ elements: bilinear square plane strain elements for the structural domain and linear elements for the acoustic domain.\
The element dimensions are $0.5\times0.5$~mm, for a total of $100 \times 100$ elements, to ensure sufficient accuracy for the acoustic and structural waves at the considered frequencies.\
The same materials and material interpolation as in~\cite{cool_exp_TO} are considered.\ 
For the structural part, this is polymethyl methacrylate (PMMA) with Young's modulus $E=4.85(1+0.05\mathrm{i})$~GPa where the imaginary part represents structural damping, density $\rho_s=1188.35$~kg/m$^3$ and Poisson's ratio $\nu=0.31$, all taken from~\cite{noe2019dynamic}.\
The core and surrounding fluid have the same properties, except for a small ratio of acoustic damping in the core through a complex speed of sound: $\rho_a=1.225$~kg/m$^3$ and $c_a=340(1+2\cdot10^{-4} \mathrm{i})$~m/s.\

\subsection*{From optimization variables to physical density fields} 
The routine uses $\bm{\xi}$ internally as design variables.\
Via a sequence of smoothing and projection filters, a physical design field is obtained.\
The smoothing filter~\cite{bourdin2001filters} is defined as follows:
\begin{linenomath}
\begin{equation}
\label{eq:dens}
\begin{aligned}
        \tilde{\xi}^e = \frac{\sum_{i=1}^{N_e} w(\bm{x}^e - \bm{x}^i)\xi^i}{\sum_{i=1}^{N_e} w(\bm{x}^e - \bm{x}^i)},\\
         w(\bm{x}^e - \bm{x}^i) = \mathrm{max}\left(0, R-||\bm{x}^e-\bm{x}^i||_2 \right),
\end{aligned}
\end{equation}
\end{linenomath}
with $\bm{x^k}=(x^k, y^k)$ the center position of element $k$.\
The smoothing filter averages the element design variables over a radius $R$.\
When this radius overlaps with the periodic boundaries at $x=0$ and $x=L_x$, the densities of elements at opposite sides of the mesh are also averaged to ensure designs are independent of the unit cell choice.\
The Heaviside projection~\cite{guest2004achieving} pushes the variables towards $0$ (fluid) and $1$ (solid).\ 
It is defined by: 
\begin{linenomath}
\begin{equation}
\label{eq:heav}
        \bar{\xi}^e = \frac{\mathrm{tanh}(\eta\beta)+\mathrm{tanh}((\xi^e-\eta)\beta)}{\mathrm{tanh}(\eta\beta)+\mathrm{tanh}((1-\eta)\beta)},
\end{equation}
\end{linenomath}
where $\eta$ is the projection level and $\beta$ is the projection strength.\\

Five distinct designs are produced: the blueprint (underscript $b$), two eroded (underscript $e$) and two dilated designs (underscript $d$).\
This isachieved using the above smoothing and projection filters, the double filtering technique~\cite{christiansen2015doublefilt} and the robust formulation~\cite{wang2011projection}, via the following operations:

\begin{linenomath}
\begin{equation}
\label{eq:double_filt}
\begin{aligned}
    &\bm{\xi} \hspace{0.2cm} \underset{R_1}{\overset{\mathrm{Eq}.~(\ref{eq:dens})}{\longrightarrow}} \hspace{0.3cm} 
    \tilde{\bm{\xi}} \hspace{0.2cm} \underset{\eta_1,\beta_1}{\overset{\mathrm{Eq}.~(\ref{eq:heav})}{\longrightarrow}} \hspace{0.2cm} 
    \bar{\tilde{\bm{\xi}}} \hspace{0.2cm} \underset{R_2}{\overset{\mathrm{Eq}.~(\ref{eq:dens})}{\longrightarrow}} \hspace{0.2cm}
    \tilde{\bar{\tilde{\bm{\xi}}}} \\ &\underset{(\eta_b,\eta_e,\eta_{e_2},\eta_d, \eta_{d_2}),\beta_2}{\overset{\mathrm{Eq}.~(\ref{eq:heav})}{\longrightarrow}} \hspace{0.2cm}
    (\bar{\tilde{\bar{\tilde{\bm{\xi}}}}}_b, \bar{\tilde{\bar{\tilde{\bm{\xi}}}}}_e, \bar{\tilde{\bar{\tilde{\bm{\xi}}}}}_{e_2}, \bar{\tilde{\bar{\tilde{\bm{\xi}}}}}_d, \bar{\tilde{\bar{\tilde{\bm{\xi}}}}}_{d_2} ),
\end{aligned}
\end{equation}
\end{linenomath}
where, for brevity, the last result is rewritten as $(\bm{\xi}_{b,\mathrm{P}}, \bm{\xi}_{e,\mathrm{P}}, \bm{\xi}_{e_2,\mathrm{P}}, \bm{\xi}_{d,\mathrm{P}}, \bm{\xi}_{d_2,\mathrm{P}})$ and $\eta_{e}=\eta_b-\Delta \eta$, $\eta_{d}=\eta_b+\Delta \eta$, $\eta_{e_2}=\eta_b-\Delta \eta_2$, $\eta_{d_2}=\eta_b+\Delta \eta_2$.\
The physical blueprint ($b$), eroded ($e$) and dilated ($d$) design fields are used for performance computation and design robustness.\
The value of $\eta$ is either kept constant at $0.1$ or varied during the optimization process until it reaches $0.1$.\
The additional eroded and dilated designs, $e_2$ and $d_2$, are used for the connectivity and volume constraints, respectively. 
$\Delta \eta_{e_2}=0.4$ and $\Delta \eta_{d_2}=0.6$ and they remain constant.\  

For the cantilever, the filter radii $R_1$ and $R_2$ are kept constant at  $10\,\mathrm{cm}$ ($20$ element widths) and $5\,\mathrm{cm}$ ($10$ element widths), respectively, and the outer frame of width $10$ elements is kept fixed at $\xi_i=1$ to achieve the outer frame.\
The filters then ensure a minimum length scale of $3.68$ element widths~\cite{QIAN2013237}.
For the sandwich panel case, the filter radii $R_1$ and $R_2$ are kept constant at $7\,\mathrm{mm}$ ($14$ element widths) and $3.5\,\mathrm{mm}$ ($7$ element widths), respectively, and the top and bottom row of 10 elements are fixed at $\xi_i=1$ to ensure a plate thickness of $5\,\mathrm{mm}$.\
Here, the filters ensure a length scale of $5.15$ element widths~\cite{QIAN2013237}.

\subsection*{STL computation}
The assumption of time-harmonic motion, weak periodicity and a unit cell size smaller than the wavelength leads to the following system of equations~\cite{fahy2007sound,boukadia2020wave}:
\begin{linenomath}
\begin{equation}
\label{eq:VA_symm}
    \left( -\omega^2 
    \underbrace{\left[ \begin{matrix} 
       \mathbf{M}_s       & 0\\
       \mathbf{S}_u         &   \mathbf{M}_a 
    \end{matrix} \right] }_{\mathbf{M}}
    +
    \underbrace{\left[ \begin{matrix} 
        \mathbf{K}_s       & \mathbf{S}_p\\
        0        &  \mathbf{K}_a 
    \end{matrix} \right] }_{\mathbf{K}} 
    \right)
    \underbrace{ \left[ \begin{matrix} 
        \mathbf{u} \\ \mathbf{p}
    \end{matrix} \right] }_{\mathbf{q}}
    =
    \underbrace{\left[ \begin{matrix} 
        \mathbf{f}_s \\ \mathbf{f}_a 
    \end{matrix} \right] }_{\mathbf{f}}
    +
    \underbrace{\left[ \begin{matrix} 
        \mathbf{e}_s \\ \mathbf{e}_a  
    \end{matrix} \right] } _{\mathbf{e}},
\end{equation}
\end{linenomath}
where $\mathbf{M}, \mathbf{K}, \mathbf{S}$ are the mass, stiffness and coupling matrices and the subscripts $s,a$ denote the structural and acoustic domain, respectively.\
Each node has three degrees-of-freedom (DOFs), namely the $x$- and $y$-displacement and the pressure. Conglomerated over all the nodes, these are denoted as the displacement vector $\mathbf{u}$ and the pressure vector $\mathbf{p}$.\
Together, they represent the DOFs vector $\mathbf{q}$, which corresponds entry-wise to the nodal force vectors $\mathbf{f}$ and $\mathbf{e}$, i.e. the internal forces and external forces, respectively.\
Infinite periodicity is enforced with Bloch-Floquet boundary conditions~\cite{bloch1929quantenmechanik}.\
Solid-fluid boundaries during optimization are dealt with via the vibro-acoustic coupling method of Jensen~\cite{jensen2019simple}.\

\subsection*{Self-weight connectivity constraint}

Connectivity is enforced with a connectivity constraint based on self-weight.\ 
For the sandwich panel, gravity points downwards, $\mu_{sw}=15$, the bottom edge ($\Gamma_b$) is fixed and periodic boundary conditions are applied on the left and right edge.\  
Furthermore, for additional restrictiveness the force is computed from the blueprint design and applied to the eroded design.
\begin{equation}
\label{eq:sw}
\left\{
\begin{array}{ll}
\hat{\mathbf{K}}_{e_2,s} \hat{\mathbf{u}}_{sw} = \hat{\mathbf{F}}_{sw} = -\mathbf{N} \bm{\xi}_{b, P} \\
\hat{\mathbf{u}}_{sw}(x) = \hat{\mathbf{u}}_{sw}(x+n L_x) \hspace{1.2cm} n = \mathbb{Z} \\
\hat{\mathbf{u}}_{sw}(x) = 0 \hspace{3.2cm} x \in \Gamma_b
\end{array}
\right. ,
\end{equation}
where $\mathbf{N} \in \mathbb{R}^{N_s/2 \times N}$ maps element design variables to structural nodal forces in the $-y$-direction and $\mathbf{u}_{sw}=\mathbf{\Lambda} \hat{\mathbf{u}}_{sw}$ is the resulting displacement vector under self-weight.\
For the cantilever, gravity points to the left,  $\mu_{sw}=4$, the left edge $ \Gamma_l$ is fixed and the gravitation force is computed on and applied to the eroded design.
\begin{equation}
\label{eq:sw_BS}
\left\{
\begin{array}{ll}
\mathbf{K}_{e_2,s} \mathbf{u}_{sw} = \mathbf{F}_{sw} = -\mathbf{N} \bm{\xi}_{e, P} \\
\mathbf{u}_{sw}(x) = 0 \hspace{3.2cm} x \in \Gamma_l
\end{array}
\right.
\end{equation}
The self-weight compliance $\theta_{sw}$ is given by:
\begin{equation}
        \theta_{sw} = \mathbf{u}_{sw}^T \mathbf{F}_{sw}.
\end{equation}
The values for $\mu_{sw}$ are selected to ensure connectivity while not overconstraining the structure towards stiffness~\cite{cool_TO_VA}.\
As the sandwich panel's performance is much more dependent on small, resonating features, the $\mu_{sw}$ value for this case is larger to accomodate this behaviour.

\subsection*{Derivatives}

The optimization problem is solved using the Method of Moving Asymptotes~\cite{svanberg1987method}.\
Since this is a gradient-based solver, derivatives for the objectives and constraints are required.\
For the STL derivatives with respect to the physical design fields, the reader is referred to appendix C1 of \cite{cool_TO_VA} for the derivation of the adjoint equations and derivatives.\ 
The self-weight compliance is self-adjoint and thus requires no additional linear system solve.\ 
A derivation of the connectivity constraint derivative can be found in appendix B.4 of \cite{cool_connectivity_review}, with the only change that the stiffness matrix depends on the blueprint but the force on the eroded design.\
The derivative of the volume constraint with respect to the physical design field $\bm{\xi}_{d_2,P}$ is straightforward.\ 
Additionally, all derivatives require the chain rule to incorporate the smoothing and projection filters and obtain the derivative with respect to $\bm{\xi}$.\ 
The derivation for this operation can be found in appendix C.2 of \cite{cool_TO_VA}.

\subsection*{Hyperparameter continuation and convergence} \label{app:continuation_convergence}

For the cantilever, the MMA control parameters are $s_{init}=0.5$, $s_{decr}=0.65$, $s_{incr}=1.07$.\
The sandwich panel uses slighty different values: $s_{init}=0.2$, $s_{decr}=0.65$, $s_{incr}=1.06$.\
Both cases use a move limit of $\mu=0.1$.\

All variations of the optimization routines in this work rely extensively on continuation, i.e., the change of hyperparameters ($\beta, ~\Delta \omega,~\ldots$) during the optimization process.\
In every variation, the hyperparameters only change after an intermediate convergence criterion is satisfied, with constant hyperparameters during every stage.\
This study considers convergence to be 1) when the constraints are satisfied (or their absolute change is less than $0.01$ between iterations) and 2) when the last $5$ iterations saw a change of less than $0.01\,\mathrm{dB}$ in the smallest STL for the sandwich panel or a relative change of the dynamic compliance less that $5.e-3$ for the cantilever.\
 When convergence is reached and there are no more hyperparameters that must be changed, the routine terminates.\
The baseline implementation uses classical $\beta$ continuation.\
For the sandwich panel, $\beta_1$ is increased from $\beta_1=1$ to, at most, $\beta_1\approx 64$, through $19$ increases by a factor of $1.2$.\
For the cantilever, $\beta_1$ is increased from $\beta_1=1$ to $2R_2=20$ in $20$ steps.
After every increase of $\beta_1$, $\beta_2$ is set to $\beta_1/2$.\ 
Note that these values are higher than the $\beta_{\mathrm{lim}}$ value proposed by da Silva et al.~\cite{daSilva2019} because, in contrast to compliance minimization, grey values are not implicitly penalized by the material interpolation scheme in dynamic topology optimization.\ 
The projection parameters $\eta_e$, $\eta_{b}$ and $\eta_{d}$ are kept constant at $0.4$, $0.5$ and $0.6$.\ 
The smoothing parameters $R_1$ and $R_2$ are also kept constant.\\

The cantilever case ends with an explicit penalization of the greyscale values.\
The Measure of Non-Discreteness $M_{nd}(\mathbf{\xi})=\sum_{i=1}^{N_e} 4\xi_i (1-\xi_i)$ for the eroded, blueprint and dilated design is added to the objective with a small factor $\lambda_{nd}$.\
That is, the new objective is  
\begin{equation}
J_{\mathrm{new}} = J_{\mathrm{old}} +\lambda_{nd} (M_{nd}(\mathbf{\xi_{e,P}}+(\mathbf{\xi_{b,P}}+(\mathbf{\xi_{d,P}})).
\end{equation}

Here, $ \lambda_{nd}$ starts at $1.e-4$ and doubles after intermediate convergence.\ 
When the change in $M_{nd}$ is less than $0.1\%$ for all three designs, the routine terminates. 

\subsection*{Updating the outer move limit}

On top of the outer move limit $\mu$ imposed internally by MMA, the bound constraints $0 \leq \bm{\xi} \leq 1$ are adapted to reflect an additional move limit $\mu_2$ such that $\Delta \bm{\xi}=\max(\bm{\xi}_{i}-\bm{\xi}_{i-1}) < \mu_2$.\ 
For the cantilever, $\mu_2=0.2$.\
For the sandwich panel optimization, $\mu_2$ is updated during the optimization due to the following two observations.\
First, when $\mu_2$ is set too low, the routine converges in a stable but slow manner.\ 
Second, when $\mu_2$ is set too high, the routine experiences numerical oscillations: STL values vary drastically between iterations.\
A value for $\mu_2$ that is as high as possible but still ensures numerical stability is desirable.\
To this end, initially and after every convergence cycle $\mu_2$ is set to $0.05$.\
STL oscillations are detected by checking if $\mathrm{STL}_b(\Delta\omega)$, $\mathrm{STL}_e(\Delta\omega)$ or $\mathrm{STL}_d(\Delta\omega)$ go from $>0.01\,\mathrm{dB}$ to $<0.01\,\mathrm{dB}$.\
If such oscillations occur for five consecutive iterations, $\mu_2$ is set to $\Delta \bm{\xi}/2$.\
On the other hand, if $\min(\mathrm{STL}_b(\Delta\omega), \mathrm{STL}_e(\Delta\omega), \mathrm{STL}_d(\Delta\omega))$ increases by more than $0.02\,\mathrm{dB}$ for five consecutive iterations, $\mu_2$ is set to $\min(3\mu_2 / 2, \mu^{\max}_2)$, where $\mu^{\max}_2$ is an adaptive upper bound set to the lowest $\mu_2$ value that caused oscillations for the current stage of hyperparameters.\\
Overall, the updating scheme allows for an automatic and cheap tuning of the move limit $\mu_2$, preventing persistent numerical oscillations while enabling fast improvements in the objective function.\
It only affects the optimization when the outer move limit is active, otherwise the internal MMA asymptotes control the convergence rate.\

\section{Monte Carlo sampling methodology} \label{app:monte_carlo_sampling_methodology}

This section elaborates on the details regarding the Monte Carlo methodology.\

\subsection*{Random initial guesses}

For the sandwich panel, the $20$ initial guesses were generated with Matlab's \verb|rand| function, which generates uncorrelated random pixel values.\ 
As this random noise is partly removed by the subsequent filtering and projection operations, the effect of a random initial guess with a spatial correlation.\
To that effect, $20$ additional guesses, now with a spatial correlation around that of the minimum length scale, were generated.\
The Monte Carlo study was repeated for half of the 16 frequency ranges.\
With this change, 3 frequency ranges changed from and to a transition region, without changing the overal conclusions of this work.\
Thus, as intended, the Monte Carlo sampling methodology removes the dependence on the initial guess.\\ 

For the cantilever, $20$ random initial guesses were generated.\
The first guess is uniform, whereas the remaining $19$ were random with a spatial correlation around the minimum length scale.\
Half of the initial guesses were enforced to be symmetric.

\subsection*{Number of optimization runs}
All optimizations are carried out on the Flemish Supercomputer (VSC). A total of $16$ frequency ranges, ranging from $[50\,\mathrm{Hz},500\,\mathrm{Hz}]$ to $[7500\,\mathrm{Hz},8000\,\mathrm{Hz}]$ are considered. All ranges have a width of $500\,\mathrm{Hz}$, except the first to prevent STL evaluation at $0\,\mathrm{Hz}$. With a sample size of $N=20$ and a total of $11$ strategies considered in this work (including the baseline), a total of $11\times20\times16=3520$ optimization runs were run.\
It is worth noting that a preliminary study using $N=10$ was performed and arrived at the same conclusions.\
This observation, together with the large computational cost, is the main limiting factor for the choice of a sample size of $N=20$.\

\subsection*{Statistical error}
The standard error of proportion ($SE$) when $N_{\text{HP}}=N_{LP}=10$ and $\tilde{P}_{\text{HP}}=50\%$ is around $\text{SE} = \sqrt{\frac{\tilde{P}_{\text{HP}}(1-\tilde{P}_{\text{HP}})}{N}}\approx 11\%$, meaning the true chance $P_{\text{HP}}$ lies within the $95\%$ confidence interval $50\% \pm 1.96*\text{SE}=[28\%,72\%]$.\ 
This wide error interval confirms that the Monte Carlo scheme is only able to distinghuish between low-performing, transition and high-performing regions.\

\clearpage
\bibliography{mybibfile}

\end{document}